\documentclass{amsart}
\usepackage{amssymb}
\usepackage{enumerate}

\def\nl{\smallskip\noindent}

\def\End{{\rm End}}
\def\Sym{{\rm S}}

\newtheorem{Thm}{Theorem}[section]
\newtheorem{Lm}[Thm]{Lemma}

\newtheorem{Def}{Definition}

\newtheorem{lemma}[Thm]{Lemma}
\newtheorem{Prop}[Thm]{Proposition}
\newtheorem{Cor}[Thm]{Corollary}

\theoremstyle{definition}

\newtheorem{Remark}[Thm]{Remark}

\newtheorem{Example}[Thm]{Example}

\renewcommand{\phi}{\varphi}
 \DeclareMathOperator{\het}{ht}

\def\Q{{\mathbb Q}}

\def\B{{\mathcal B}}

\newcommand{\eps}{\varepsilon}

\newcommand{\A}{{\rm A}}
\newcommand{\D}{{\rm D}}
\newcommand{\E}{{\rm E}}

\newcommand{\adj}{\sim}

\def\a{\alpha}
\def\b{\beta}
\def\c{\gamma}

\def\gam{\gamma}
\def\alp{\alpha}

\def\t{\tau}

\author{Arjeh M. Cohen
\& Di\'{e} A.H. Gijsbers \& David B. Wales}
\address{Arjeh M. Cohen\\
Department of Mathematics and Computer Science\\
Eindhoven University of Technology\\
POBox 513\\
5600 MB Eindhoven\\
The Netherlands} \email{A.M.Cohen@tue.nl}
\address{Di\'{e} A.H. Gijsbers\\
Department of Mathematics and Computer Science\\
Eindhoven University of Technology\\
POBox 513\\
5600 MB Eindhoven\\
The Netherlands} \email{D.A.H.Gijsbers@tue.nl}
\address{David B. Wales\\
Department of Mathematics\\
Caltech\\
Pasadena, CA 91125\\
USA} \email{dbw@its.caltech.edu}

\title{A Poset Connected to Artin Monoids of Simply Laced Type}

\begin{document}

\begin{abstract}
Let $W$ be a Weyl group whose type is a simply laced Dynkin
diagram. On several $W$-orbits of sets of mutually commuting
reflections, a poset is described which plays a role in linear
representations of the corresponding Artin group $A$. The poset
generalizes many properties of the usual order on positive roots
of $W$ given by height. In this paper, a linear representation of
the positive monoid of $A$ is defined by use of the poset.
\end{abstract}

\maketitle

\section{Introduction}\label{sec:intro}
The beautiful properties of the high root used in \cite{CGW} to
construct Lawrence-Krammer representations of the Artin group with
non-commutative coefficients have analogues for certain sets of
orthogonal roots.  We study these properties and exploit them to
construct a linear representation of the Artin monoid. In many
instances, the monoid representation extends to an Artin group
representation; this will be the subject of subsequent work.

Let $M$ be a Coxeter diagram of simply laced type, i.e., its
connected components are of type $\A$, $\D$ or $\E$. The
Lawrence-Krammer representation (\cite{bigelow,CW,digne,Krammer})
has a basis consisting of positive roots of the root system of the
Weyl group $W=W(M)$ of type $M$.  Here we use instead, a $W$-orbit
$\B$ of sets of mutually orthogonal positive roots. Not all
$W$-orbits of this kind are allowed; we call those which are
allowed, {\em admissible} (cf.\ Definition \ref{def:admissible}; a
precise list for $M$ connected is in Table~\ref{tab:type}). In the
Lawrence-Krammer representation we used the partial ordering of
the positive roots given by $\b\le \c$ iff $\c-\b$ is a sum of
positive roots with non-negative coefficients. In Proposition
\ref{Prop:AdmOrd} we generalize this ordering to an ordering
$(\B,<)$ on admissible $W$-orbits $\B$ of mutually orthogonal
roots. In the action of $w\in W$ on a set $B\in\B$ the image $wB$
is the set of positive roots in $\{\pm w\b\mid \b\in B\}$.  In the
single root case, there is a unique highest element, the
well-known highest root. This property extends to $(\B,<)$: there
is a unique maximal element $B_0$ in $\B$ (cf.\ Corollary
\ref{uniquemax}).

In the Lawrence-Krammer representation, the coefficients were
obtained from the Hecke algebra whose type is the subdiagram of
$M$ induced on the set of nodes $i$ of $M$ whose corresponding
fundamental root $\a_i$ is orthogonal to the highest root.  Here,
the coefficients are obtained from the Hecke algebra $Z$ whose
type is the  subdiagram of $M$ induced on the nodes $i$ whose
corresponding fundamental root $\a_i$ is orthogonal to each
element of $B_0$. Moreover, in the Lawrence-Krammer
representation, to each pair of a positive root $\b$ and a node
$i$ with corresponding fundamental root $\a_i$ such that
$(\a_i,\b) = 0$, we assigned an element $h_{\b,i}$ of the
coefficient algebra. It occurs in the definition of the action of
a fundamental generator of the Artin group $A$ in the
Lawrence-Krammer representation, on the basis element
 $\b$. For the analogous purpose, we introduce
elements $h_{B,i}$ (Definition \ref{defw}) in the corresponding
coefficient algebra $Z$. These elements are parameterized by pairs
consisting of an element $B$ of $\B$ and a node $i$ of $M$ such
that the corresponding fundamental root $\a_i$ is orthogonal to
all of $B$.

In analogy to the developments in \cite{CGW} we define a right
free $Z$-module with basis $x_B$ ($B\in \B$) which is a left
module for the positive monoid $A^+$ of the Artin group $A$ of
type $M$.  For each node $i$ of $M$, the $i$-th fundamental
generator $s_i$ of $A^+$ maps onto the linear transformation
$\tau_i$ on $V$ given by the following case division.
\begin{eqnarray}\label{taui}
\tau_ix_B &=& \left\{
\begin{tabular}{rrlc}
$0$  && if $\a_i \in B $,   &\\
$x_Bh_{B,i}$  && if $\a_i \in B^{\perp}$,    &\\
$x_{r_iB}$ && if $r_iB<B$,        &\\
$x_{r_iB}-mx_B$ && if $r_iB>B$.     &
\end{tabular}
\right.
\end{eqnarray}
This leads to the following main result of this paper.

\begin{Thm}\label{th:main}
Let $W$ be a Weyl group of simply laced type. For $\B$ an
admissible $W$-orbit of sets of mutually orthogonal positive
roots, there is a partial order $<$ on $\B$ such that the above
defined map $s_i\mapsto \tau_i$ determines a homomorphism of
monoids from $A^+$ to $\End(V)$.
\end{Thm}

In the sections below, we deal with this construction in detail,
the proof of the theorem is in Section 5.

When labeling the nodes of an irreducible diagram $M$, we will
choose the labeling of \cite{Bourb}. If $M$ is disconnected, the
representations are easily seen to be a direct sum of
representations corresponding to the components.  Since the poset
construction also behaves nicely, it suffices to prove the theorem
only for $M$ connected.  Therefore, we will assume $M$ to be
connected for the greater part of this paper.

\section{Admissible orbits}
Let $M$ be a spherical Coxeter diagram.  Let $(W,R)$ be the
Coxeter system of type $M$ with $R= \{r_1,\ldots,r_n\}$.
Throughout this paper we shall assume that $M$ is simply laced,
which means that the order of each product $r_ir_j$ is at most
$3$.

By $\Phi^+$ we denote the positive root system of type $M$ and by
$\a_i$ the fundamental root corresponding to the node $i$ of $M$.
We are interested in sets $B$ of mutually commuting reflections.
Since each reflection is uniquely determined by a positive root,
the set $B$ corresponds bijectively to a set of mutually
orthogonal roots of $\Phi^+$.  We will almost always identify $B$
with this subset of $\Phi^+$.  The action of $w\in W$ on $B$ is
given by conjugation in case $B$ is described by reflections and
by $w\{\b_1,\ldots,\b_p\} = \Phi^+\cap\{\pm w\b_1,\ldots,\pm
w\b_p\}$ in case $B$ is described by positive roots.  The action
of an element $w\in W$ on $B$ should not be confused with the
action of $w$ on a root: in our case we have $w\{\a_i\} =
\{\a_i\}$ whereas the usual action on roots implies $w\a_i =
\a_i$.  For example, if $r_i$ is the reflection about $\a_i$,
$r_i\{\a_i\}=\{\a_i\}$ but $r_i\a_i=-\a_i$.

The $W$-orbit $\B$ of a set $B$ of mutually orthogonal positive
roots is the vertex set of a graph with edges labeled by the nodes
of $M$, the edges with label $j$ being the unordered pairs
$\{B,r_jB\}$ (so $r_jB\ne B$) for $B\in \B$. The results of
Section 3 show that if $\B$ is admissible, the edges of this graph
can be directed so as to obtain a partially ordered set (poset)
having a unique maximal element. This section deals with the
notion of admissibility.

We let $\het(\b)$ be the usual height of a root $\b \in \Phi^+$
which is $\sum_i a_i$ where $\b=\sum a_i\a_i$.

\begin{Prop}\label{DnsystemIsAdm}
Let $M$ be of simply laced type. Every $W$-orbit of sets of
mutually orthogonal positive roots satisfies the following
properties for $B\in \B$, $j\in M$ and $\b,\c\in B$.
\begin{enumerate}[(i)]
\item There is no node $i$ for which $(\a_i,\b)=1,$ $(\a_i,\c)=-1$
and $\het(\b)=\het(\c)+1$. \item Suppose $(\a_j,\b)=-1$ and
$(\a_j,\c)=1$ with $\het(\c)=\het(\b)+2$. Then there is no node
$i$ for which $\a_i\in B^\perp$ and $i\sim j$.
\end{enumerate}
\end{Prop}

\begin{proof}
Let $B$ be a set of mutually orthogonal positive roots, and
$\b,\gamma\in B$.

\nl(i). Suppose there is a node $i$ for which $(\a_i,\b)=1,$
$(\a_i,\c)=-1$ and $\het(\b)=\het(\c)+1$. As $\b$ and $\c$ are
orthogonal we have $(\b,\c+\a_i)= 1$ so $\b-\gamma-\a_i \in \Phi$.
This is not possible as $\het(\b - \c- \a_i) = 0$.

\nl(ii). Let $\b$ and $\c$ be as in the hypothesis and assume
there is an $i$ for which $\a_i\in B^\perp$ and $i\sim j$. Then
$(\a_i,\c-\a_j) = 1$, so $\c-\a_j-\a_i$ is a root. As
$\het(\c)=\het(\b)+2$ we have $\het(\c-\a_j-\a_i)=\het(\b)$. But
$(\b,\c-\a_j-\a_i) = 1$, so $\b-\c+\a_j+\a_i$ is a root which
contradicts $\het(\b-\c+\a_j+\a_i)=0$.
\end{proof}

\begin{Def}\label{def:admissible}
Let $\B$ be a $W$-orbit of sets of mutually orthogonal positive
roots. We say that $\B$ is {\em admissible} if for each $B\in \B$
and $i,j\in M$ with $i\not\sim j$ and $\c,\c-\a_i+\a_j \in B$, we
have $r_iB=r_jB$.
\end{Def}

Not all $W$-orbits on sets of mutually orthogonal positive roots
are admissible. The $W$-orbit of the triple
$\{\a_1+\a_2+\a_3,\a_2+\a_3+\a_4,\a_1+\cdots+\a_5\}$ of positive
roots for $M = \D_5$ is a counterexample.  Suppose that $M$ is
disconnected with components $M_i$.  Then $\B$ is admissible if
and only if each of the corresponding $W(M_i)$-orbits is
admissible. So there is no harm in restricting our admissibility
study to the case where $M$ is connected. In that case,
Proposition~\ref{prop:adm-examples} below gives a full
characterization of admissible orbits.

\begin{Lm}\label{lm:otherreflection}
Let $r$ be a reflection in $W$ and let $\b,\gam$ be two mutually
orthogonal positive roots moved by $r$. Then there exists a
reflection $s$ which commutes with $r$ such that
$\{\b\}=rs\{\gam\}$.
\end{Lm}

\begin{proof}
Let $\delta_r$ be the positive root corresponding to the
reflection $r$. Now $r\{\b\} = \{\pm \b \pm \delta_r\}\cap \Phi^+$
and $(\gam, \b\pm \delta_r) = \pm(\gam, \delta_r)$. Using this we can construct
a new positive root $\delta$ depending on
$(\b,\delta_r)$, $(\gam,\delta_r)$
as indicated in the table below.

\begin{center}
\begin{tabular}{|c|c|c|c|}
\hline
$(\b,\delta_r)$ & $(\gam,\delta_r)$ & $\pm \delta$ \\
\hline
$1$ & $1$   & $\b+\gam-\delta_r$ \\
$1$ & $-1$  & $\b-\gam-\delta_r$ \\
$-1$ & $1$  & $\b-\gam+\delta_r$ \\
$-1$ & $-1$ & $\b+\gam+\delta_r$\\
\hline
\end{tabular}
\end{center}

It is easy to check that the reflection $s$ with
root $\delta$ commutes with $r$ and indeed $\{\b\}= rs\{\gam\}$.
\end{proof}

\begin{Prop}\label{prop:adm-examples}
Let $M$ be connected. The following statements concerning a
$W$-orbit $\B$ of sets of mutually orthogonal positive roots are
equivalent.

\begin{enumerate}[(i)]
\item $\B$ is admissible. \

\item For each pair $r$, $s$ of commuting reflections of $W$ and
each $B\in \B$ such that $\gamma$ and $rs\gamma$ both belong to
$B$, we have $rB = sB$.

\item For each reflection $r$ of $W$ and each $B\in \B$ the size
of $rB\setminus B$ is one of $\ 0,1,2,4$.

\item $\B$ is one of the orbits listed in Table~\ref{tab:type}.
\end{enumerate}
\end{Prop}

\noindent Below in the proof we show that four is the maximum
possible roots in $rB\setminus B$ which can be moved and so only
three is ruled out in part (iii).

\begin{proof}
(i)$\implies$(ii). By (i), assertion (ii) holds when $r$ and $s$
are fundamental reflections. The other cases follow by conjugation
since each pair of commuting reflections is conjugate to a pair of
fundamental conjugating reflections. (As each reflection is
conjugate to a fundamental reflection, the reflections orthogonal
to it can be determined and the system of roots orthogonal to a
reflection has the type obtained by removing nodes connected to
the extending node of the affine diagram.)

\medskip

\nl (ii)$\implies$(iii). When all $r,s \in W$ move at most two
mutually orthogonal roots, the implication holds trivially. If $r$
would move five mutually orthogonal roots then the $6\times6$ Gram
matrix for these roots together with the root of $r$ is not
positive semi-definite as its determinant is $-16$, a contradiction. Hence $r$ moves at most $4$ roots.

Assume we have a $B \in \B$ such that $r$ moves precisely three roots
of $B$, say $\b_1,\b_2,\b_3$. By Lemma~\ref{lm:otherreflection} we
know there exists a reflection $s$ such that $\b_1=rs\b_2$. Now $\b_2
= \b_1 \pm \delta_r \pm \delta_s$ with $\delta_r$, $\delta_s$ the
positive roots corresponding to $r$ and $s$, respectively.
As $\b_3$ is orthogonal to $\b_1$ and $\b_2$, we find
$(\b_3,\delta_s) = \pm (\b_3,\delta_r)$, so $s$ moves $\b_3$ as
well. But obviously $r\b_3 \ne s\b_3$, so $rB\ne sB$, which contradicts
(ii).

\nl (iii)$\implies$(i). Let $B\in \B$ and $i,j\in M$ with
$i\not\sim j$ and $\c,\c-\a_i+\a_j \in B$. When both $r_i,r_j$ do
not move any other root then $r_iB=r_jB$. Without loss of
generality we can assume $r_i$ moves four roots of $B$. Let $\b$
be a third root in $B$ moved by $r_i$. As $\b$ has to be
orthogonal to $\c,\c-\a_i+\a_j$ we find $(\a_i,\b)=(\a_j,\b)$, so
$\b-\a_i-\a_j$ or $\b+\a_i+\a_j$ is a positive root as well. This
root is also moved by $r_i$ and mutually orthogonal to
$\c,\c-\a_i+\a_j$ and $\b$.

So now $\{\c,\c-\a_i+\a_j, \b,\b-\a_i-\a_j\} \subseteq B$. But these
$4$ roots are also the roots moved by $r_j$.  We know from above that $4$
is the maximal number of mutually orthogonal roots moved by $r_i$ (or by $r_j$ for
that matter).  We find $r_iB=r_jB$ which proves $\B$ is
admissible.

At this point we have achieved equivalence of (i), (ii), and (iii), a
fact we will use throughout the remainder of the proof.

\nl (iii)$\implies$(iv). In Table~\ref{tab:orbits}, we have listed
all $W$-orbits of sets of mutually orthogonal positive roots. It
is straightforward to check this (for instance by induction on the
size $t$ of such a set), so we omit the details. For all orbits in
Table~\ref{tab:orbits} but not in Table~\ref{tab:type} we find,
for some set $B$ in the orbit $\B$, a reflection $r$ which moves
precisely three roots.

\begin{center}
\begin{table} \begin{tabular}{|c|c|c|}
\hline
$M$  & $|B|$ & $B$ \\
\hline
$\A_n$ & $t$  & $\{\a_1,\a_{3},\ldots,\a_{2t-1}\}$ \\
\hline $\D_n$ & $t$  & $\{\a_i,\b_{i} \mid i=1,3,\ldots,..,2k-1\}
\quad \cup \quad \{\a_i \mid i=2k+1,2k+3, \ldots,..,2t-1 \}$  \\
$\D_n$ & $n/2$  & $\{\a_1,\a_3,\ldots,\a_{n-3},\a_{n}\}$  \\
\hline
$\E_6$ & $1$  & $\{\a_2\}$    \\
$\E_6$ & $2$  & $\{\a_2,\a_5\}$   \\
$\E_6$ & $3$  & $\{\a_2,\a_3,\a_5\}$   \\
$\E_6$ & $4$  & $\{\a_2,\a_3,\a_5,\a_2+\a_3+2\a_4+\a_5\}$  \\
\hline
$\E_7$ & $1$  & $\{\a_2\}$ \\
$\E_7$ & $2$  & $\{\a_2,\a_5\}$   \\
$\E_7$ & $3$  & $\{\a_2,\a_5,\a_7\}$   \\
$\E_7$ & $3$  & $\{\a_2,\a_3,\a_5\}$   \\
$\E_7$ & $4$  & $\{\a_2,\a_3,\a_5,\a_7\}$  \\
$\E_7$ & $4$  & $\{\a_2,\a_3,\a_5,\a_2+\a_3+2\a_4+\a_5\}$  \\
$\E_7$ & $5$  & $\{\a_2,\a_3,\a_5,\a_7,\a_0\}$  \\
$\E_7$ & $6$  & $\{\a_2,\a_3,\a_5,\a_7,\a_2+\a_3+2\a_4+\a_5,\a_0\}$  \\
$\E_7$ & $7$  & $\{\a_2,\a_3,\a_5,\a_7,\a_2+\a_3+2\a_4+\a_5,\a_2+\a_3+2\a_4+2\a_5+2\a_6+\a_7,\a_0\}$ \\
\hline
$\E_8$ & $1$  & $\{\a_2\}$ \\
$\E_8$ & $2$  & $\{\a_2,\a_5\}$   \\
$\E_8$ & $3$  & $\{\a_2,\a_3,\a_5\}$   \\
$\E_8$ & $4$  & $\{\a_2,\a_3,\a_5,\a_7\}$   \\
$\E_8$ & $4$  & $\{\a_2,\a_3,\a_5,\a_2+\a_3+2\a_4+\a_5\}$   \\
$\E_8$ & $5$  & $\{\a_2,\a_3,\a_5,\a_7,\a_0\}$   \\
$\E_8$ & $6$  & $\{\a_2,\a_3,\a_5,\a_7,\a_0,\bar{\a}_0\}$   \\
$\E_8$ & $7$  & $\{\a_2,\a_3,\a_5,\a_7,\a_2+\a_3+2\a_4+\a_5,\a_0,\bar{\a}_0\}$   \\
$\E_8$ & $8$  & $\{\a_2,\a_3,\a_5,\a_7,\a_2+\a_3+2\a_4+\a_5,\a_2+\a_3+2\a_4+2\a_5+2\a_6+\a_7,\a_0,\bar{\a}_0\}$ \\
\hline
\end{tabular}
\bigskip
\caption{For $M=\A_n$ we have $t \le \frac{n+1}{2}$ and for
$M=\D_n$ we have $t \le \frac{n}{2}$. For $M=\D_n$ we write
$\b_{n-1} = \a_n$ and $\b_{2t+1} = \a_n+\a_{n-1}+2\a_{n-2}+ \ldots
+ 2\a_{2t+2} + \a_{2t+1}$. In $E_7$ and $E_8$ we use $\a_0 =
2\a_1+2\a_2+3\a_3+4\a_4+3\a_5+2\a_6+\a_7$ and $\bar{\a}_0 =
2\a_1+3\a_2+4\a_3+6\a_4+5\a_5+4\a_6+3\a_7+2\a_8$, the respective
highest roots.}  \label{tab:orbits}

\end{table}
\end{center}

\bigskip
We will use the observation that if $B$ belongs to a
non-admissible orbit for $W$ of type $M$, then it also does not
belong to an admissible orbit for $W$ of any larger type.

\nl For $M = \D_n$ the sets not in Table~\ref{tab:type} contain at
least one pair of roots $\eps_i-\eps_j,\eps_i+\eps_j$ but also at
least one root $\eps_p \pm \eps_q$ without the corresponding other
positive root containing $\eps_p,\eps_q$. (Here the $\eps_i$ are
the usual orthogonal basis such that $\Phi^+ =\{ \eps_i \pm \eps_j
\mid i < j\}$.) For these sets the reflection corresponding to a
positive root $\eps_j \pm \eps_q$ moves precisely three roots.

\nl Suppose $M= \E_n$. The orbit of sets of three mutually orthogonal
roots which is not in Table~\ref{tab:type} is the orbit of
$\{\a_2,\a_3,\a_5\}$, which is not admissible in the subsystem of
type $\D_4$ corresponding to these three roots and $\a_4$, as
$r_4$ moves all three roots.

The orbit of four mutually orthogonal positive roots not in
Table~\ref{tab:type} contains the set $\{\a_2,\a_3,\a_5,\a_7\}$
and $r_4$ moves exactly three of these.

The orbit of five mutually orthogonal positive roots not in
Table~\ref{tab:type} contains the set
$\{\a_2,\a_3,\a_5,\a_7,\a_0\}$ and $r_4$ moves again exactly three
of these.

If $M=\E_7$, the orbit of sets of six mutually orthogonal positive
roots containing
$\{\a_2,\a_5,\a_7,\a_3,\a_2+\a_3+2\a_4+\a_5,\a_0\}$ remains.
Clearly the reflection $r_1$ moves only the last three roots. If
$M = \E_8$, the orbit of sets of six mutually orthogonal positive
roots is not admissible as it contains the orbit of $\E_7$ we just
discussed.

Finally the orbit of seven mutually orthogonal positive roots in
$\E_8$ contains
$\{\a_2,\a_3,\a_5,\a_2+\a_3+2\a_4+\a_5,\a_7,\a_0,\bar{\a}_0\}$.
Here the reflection $r_8$ moves only the last three roots.

\nl (iv)$\implies$(iii). All orbits for type $\A_n$ are admissible
as here every reflection moves at most two mutually orthogonal
roots.

All sets in the first collection of orbits in $\D_n$ contain from
every pair of roots $\eps_i-\eps_j,\eps_i+\eps_j$ at most one
element. So again as for $\A_n$, every reflection moves at
most two mutually orthogonal roots.

All sets in the second collection of orbits in $\D_n$ contain from
every pair of roots $\eps_i-\eps_j,\eps_i+\eps_j$ both roots or
none of them. So every reflection here will always move an even
number of roots, so the size of $rB\setminus B$ is equal to $0,2$
or $4$.

The orbits in $\E_n$ of sets containing fewer than three roots are
admissible as every reflection will never move more than two
roots. For the remaining six orbits it is an easy exercise to
verify for one chosen set that every reflection moves indeed
$0,1,2$ or $4$ roots.
\end{proof}

\bigskip

\begin{table}\begin{tabular}{|c|c|c|c|c|c|}
\hline
$M$  & $|B|$ & $|\B|$ & $Y$ & $C$ & $N_W(B)$\\
\hline
$\A_n$ & $t$  & $\frac{(n+1)!}{2^tt!(n-2t+1)!}$  & $\A_{n-2t}$ & $\A_{n-2t}$ & $2^t\Sym_{t}\Sym_{n+1-2t}$\\
\hline
$\D_n$ & $t$& $\frac{n!}{t!(n-2t)!}$ & $\A_1^t\D_{n-2t}$ & $\A_1\D_{n-2t}$ & $2^{2t}\Sym_{t}W(\D_{n-2t})$ \\
$\D_n$ & $2t$& $\frac{n!}{2^{t}t!(n-2t)!}$ & $\D_{n-2t}$ & $\A_{n-2t-1}$ & $2^{2t}W({\rm B}_{t})W(\D_{n-2t})$\\
\hline
$\E_6$ & $1$  & $36$   & $\A_5$ & $\A_5$ & $2\Sym_6$\\
$\E_6$ & $2$  & $270$  & $\A_3$ & $\A_2$ & $2^{2+1}\Sym_4$ \\
$\E_6$ & $4$  & $135$  & $\emptyset$ & $\emptyset$ & $2^4\Sym_4$ \\
\hline
$\E_7$ & $1$  & $63$  & $\D_6$ & $\D_6$ & $2W(\D_6)$ \\
$\E_7$ & $2$  & $945$ & $\A_1\D_4$ & $\A_1\D_4$ & $2^{2+1+1}W(\D_4)$ \\
$\E_7$ & $3$ & $315$  & $\D_4$ & $\A_2$ & $2^{3}\Sym_3W(\D_4)$ \\
$\E_7$ & $4$  & $945$  & $\A_1^3$ & $\A_1$ & $2^{4+3}\Sym_4$ \\
$\E_7$ & $7$  & $135$  & $\emptyset$ & $\emptyset$ & $2^{7}{\rm L}(3,2)$ \\
\hline
$\E_8$ & $1$  & $120$  & $\E_7$ & $\E_7$ & $2W(\E_7)$ \\
$\E_8$ & $2$  & $3780$ & $\D_6$ & $\A_5$ & $2^{2+1}W(\D_6)$ \\
$\E_8$ & $4$  & $9450$ & $\D_4$ & $\A_2$ & $2^4\Sym_3W(\D_4)$ \\
$\E_8$ & $8$  & $2025$  & $\emptyset$ & $\emptyset$ & $2^{8+3}{\rm L}(3,2)$ \\
\hline
\end{tabular}

\bigskip
\caption{Each row contains the type $M$, the size of $B \in \B$,
the size of the $W$-orbit $\B$ containing $B$, the Coxeter type
$Y$ of the roots orthogonal to $B$, the type of the Hecke Algebra
$C$ defined in Corollary \ref{cor:hB}, and the structure of the
normalizer in $W$ of $B$, respectively. In the first line for
$\D_n$, we define $\D_{n-2t}$ as being empty if $n-2t\le1$.  Only one of the 
roots $\eps_i\pm \eps_j$ occur for roots in the first line of $D_n$.  For roots 
in the second line, both occur. }
\label{tab:type}
\end{table}
\bigskip

We finish this section with some further comments on
Table~\ref{tab:type}. If $B$ is a set of mutually orthogonal
positive roots as indicated in Table~\ref{tab:type}, then the type
$Y$ of the system of all roots orthogonal to $B$ is listed in the
table. In the final column of the table we list the structure of
the stabilizer in $W$ acting on $\B$. If $\B$ has an element $B$
all of whose members are fundamental roots, this stabilizer can be
found in \cite{Howlett}. Two distinct lines represent different
classes; sometimes even more than two, in which case they fuse
under an outer automorphism (so they behave identically). This
happens for instance for $M=\D_{n}$ (first line) with $n=2t$. In
the second line for $\D_n$, the permutation action is not
faithful.

\section{Posets}\label{sec:posets}
In this section we show that admissible orbits carry a nice poset structure.
An arbitrary $W$-orbit of sets of mutually orthogonal positive roots satisfies
all of the properties of the proposition below except for (iii).

\begin{Prop}\label{Prop:AdmOrd}
Let $M$ be a spherical simply laced diagram and $\B$ an admissible
$W$-orbit of sets of mutually orthogonal positive roots. Then there is an
ordering $<$ on $\B$ with the following properties.
\begin{enumerate}[(i)]
\item For each node $i$ of $M$ and each $B\in \B$, the sets $B$
and $r_i B$ are comparable. Furthermore, if $(\a_i,\b)=\pm1$ for some
$\b \in B$, $r_iB\ne B$.
\item Suppose $i\sim j$ and $\a_i\in
B^\perp$. If $r_jB < B$, then $r_ir_jB<r_jB$. Also, $r_jB > B$
implies $r_ir_jB > r_j B$. \item If $i\not\sim j$, $r_i B<B$,
$r_jB < B$, and $r_i B\ne r_j B$, then $r_ir_jB<r_jB$ and
$r_ir_jB<r_i B$. \item If $i\sim j$, $r_i B<B$, and $r_jB < B$,
then either $r_ir_jB = r_jB$ or $r_ir_jB<r_jB$, $r_jr_iB<r_i B$,
$r_ir_jr_iB< r_ir_jB$, and $r_ir_jr_iB<r_jr_iB$.
\end{enumerate}
\end{Prop}

It readily follows from the existence result that there is a
unique minimal ordering $<$ satisfying the requirements of the
proposition (it is the transitive closure of the pairs $(B,r_jB)$
for $B\in \B$ and $j$ a node of $M$ such that $r_jB>B$).  This
poset $(\B,<)$ with this minimal ordering is called the {\em
monoidal poset} (with respect to $W$) on $\B$ (so $\B$ should be
admissible for the poset to be monoidal).

\begin{proof}
We define the relation $<$ on $\B$ as follows: for $B,C\in \B$ we
have $B<C$ iff there are $\b\in B\setminus C$ and $\c\in
C\setminus B$, of minimal height in $B\setminus C$, respectively
$C\setminus B$, such that $\het(\b)<\het(\c)$. It is readily
verified that $<$ is an ordering. We show that it also satisfies
properties (i),...,(iv).

We will need various properties involving the actions of the $r_i$
on $\B$. Clearly, if $\a_i\in B^\perp$, then $r_iB = B$.  As
described earlier, if $\a_i$ is in $B$ we replace $r_i\a_i=-\a_i$
with $\a_i$ and say $\a_i$ is fixed.  Note that then also $r_iB =
B$.  If $(\a_i,\b)=\pm1$, we say $r_i$ moves $\b $.   In this case with
Proposition \ref{Prop:AdmOrd}(i) in mind, we see $r_iB\neq B$ and we say $r_i$ {\em
lowers} $B$ if $B>r_iB$ and we say $r_i$ {\em raises} $B$ if
$B<r_iB$.  We also use this for a single root $\b$: if $\b+\a_i$
is a root, we say $r_i$ raises $\b$ and if $\b-\a_i$ is a root we
say $r_i$ lowers $\b$.

\nl(i).  If $\a_i$ is orthogonal to each member of $B$ or $\a_i\in
B$, then $r_i B = B$, so $B$ and $r_i B$ are comparable.  So we
may assume that $(\a_i,\b) =\pm 1$ for at least one $\b$.  Notice
if $(\b,\a_i)=\pm1$, that $r_i\b = \b \pm \a_i$ is in $\Phi$ but
not in $B$ as $(\b,\b \pm \a_i)=2\pm 1\ne 0$ whereas different
elements of $B$ are orthogonal.  In particular, $r_iB\neq B$.  If $(\b,\a_i)=\pm 1$ holds for
exactly one member of $B$, then clearly $r_i B$ and $B$ are
comparable. Suppose now $B$ and $r_i B$ are not comparable. Then
there exist at least one $\b \in B\setminus r_iB$  and one $\c\in
B$ such that $\het(\b) = \het(r_i\c)$ and $r_i\c$ is an element of
minimal height in $r_iB\setminus B$.  Clearly $\b \ne r_i\c$ as
they are in different sets by their definition. As $\c \in
B\setminus r_iB$ we have $\het(\c) \ge \het(\b)$ by our assumption
that $\b$ is of minimal height in $B\setminus r_iB$.  As $\het
r_i\c = \het \b$ we have $(\a_i,\c)=1$. Also $r_i\b \in
r_iB\setminus B$ and as $r_i\c$ with $\het(r_i\c) = \het(\b)$ is
of minimal height in $r_iB\setminus B$ we see $\het(r_i\b)\ge
\het(r_i\c)$.  In particular we must have $(\a_i,\b)=-1$. But
according to Condition (i) Proposition \ref{DnsystemIsAdm}, this never occurs.

\nl(ii). By the assumption $r_jB<B$, there is a root $\b \in B$ of
minimal height among those moved by $r_j$ such that $\b-\a_j\in
r_jB$. Then $\b-\a_j$ is minimal among those moved by $r_i$ in
$r_jB$ and $\b-\a_i-\a_j\in r_ir_jB$, so $r_ir_jB<r_jB$.

The proof for the second assertion is a bit more complicated.  By
the assumption $r_jB>B$, there is no root $\b \in B$ of minimal
height among those moved by $r_j$ such that $\b-\a_j\in r_jB$.
Indeed all that are moved go to $\b+\a_j$ in $r_jB$.  Suppose that
$\delta$ has minimal height among the roots moved by $r_i$ in $r_j
B$.  This implies $\delta=\c\pm \a_j$ for some $\c\in B$. If
$\delta=\c+\a_j$ for all choices of $\delta$, then $r_ir_jB>r_jB$,
as $r_i(\c+\a_j)=\c+\a_j+\a_i$.  So assume that $\delta=\c-\a_j$
has minimal height among the roots moved by $r_i$ in $r_jB$ for
some $\c\in B$.  Let $h$ be the minimal height of all elements of
$B$ moved by $r_j$.  We know each of these roots is raised in
height by $r_j$ and so $\c$ is not one of them.  In particular
$\het(\c)>h$.  Also $\het(\c)-1 = \het(\delta) \le h+1$.  It
follows that $\het(\c) \le h+2$. The two cases are $\c$ has height
$h+1$ or $h+2$.  By Condition (i) for Proposition \ref{DnsystemIsAdm} for $\c $ and
$\b$ the case $h+1$ is ruled out.  But then Condition (ii) of
Proposition \ref{DnsystemIsAdm} with $\a_i$, $\c$, and $\b$ rules out the case
$h+2$.

\nl(iii).  Suppose $r_jB<B$ and $r_iB<B$ with $r_iB\ne r_jB$.
Choose $\b$ an element of smallest height in $B$ moved by $r_j$
for which $\b-\a_j$ is a root.  Choose $\c$ an element of smallest
height in $B$ moved by $r_i$ with $\c-\a_i$ a root.  We are
assuming $\b-\a_j$ is a root.  This means as $(\a_i,\a_j)=0$ that
$\b\pm \a_i$ is a root if and only if $\b-\a_j\pm \a_i$ is a root.

To prove the result we will get a contradiction if we assume $r_i$
raises $r_jB$.  Suppose then $r_i$ raises $r_jB$.  In this case
all elements $\zeta$ of smallest height in $r_jB$ which are moved
by $r_i$ have $\zeta+\a_i$ as roots.  We will show first that
$\c-\a_j$ is not a root.  If it were, $r_i$ lowers it as $\c-\a_i$
is a root.  This means it is a root of smallest height moved by
$r_i$ as $\c$ is a root of smallest height moved by $r_i$ in $B$
and in $r_jB$ this has height one smaller.  But it is lowered, not
raised.  This means $\c-\a_j$ is not a root.

Depending on $(\c,\a_j)$, either $\c$ or $\c+\a_j$ is a root of
$r_jB$. Suppose $(\c,\a_j)=0$ and so $\c$ is a root of $r_jB$.  As
$r_i$ raises $r_jB$, all elements of smallest height moved by
$r_i$ must be raised.  As $\c$ is lowered, there must be an
$r_j\delta\in r_jB$ with $\delta-\a_j$ a root of $r_jB$ and
$\het(\delta-\a_j)$ less than $\het(\c)$.  Its height must be one
less than $\het(\c)$ as heights are lowered at most one by $r_j$.
Now in $r_jB$, the elements $\delta-\a_j$ and $\c$ contradict
condition (i) for Proposition \ref{DnsystemIsAdm}.  Suppose then
$\c+\a_j$ is root. The smallest height of elements for which $r_i$
moves roots in $r_jB$ is now either $\het (\c)$ or $\het (\c)-1$.
(It cannot be $\het(\c)+1$ as $\c+\a_j$ is lowered.)  If it is
height $\het(\c)$ there is an element $\delta$ of height
$\het(\c)$ which is raised by $r_i$.  Now $\delta$ and $\c+\a_j$
contradict Condition (i) of Proposition \ref{DnsystemIsAdm}.

We are left with one case in which an element of height
$\het(\c)-1$ in $r_jB$ is raised by $r_i$.  This means there is
$\delta $ in $B$ of height $\het(\c)$ which is lowered by $r_j$
and raised by $r_i$. Recall $\c$ is lowered by $r_i$ and raised by
$r_j$.
As $i \not \adj j$ we have $\delta+\a_i-\a_j, \c+\a_j-\a_i$ in
$r_ir_jB$.  Now $(\delta, \c+\a_j-\a_i) = 2$ so $\delta = \c
+\a_i-\a_j$.
By admissibility of $\B$ we have $r_iB=r_jB$ contradicting the
starting assumptions.

\nl(iv). We shall use the results of the following computations
throughout the proof. Let $\epsilon \in B$ and write $\rho =
(\a_i,\epsilon)$ and $\sigma = (\a_j,\epsilon)$. Then
$\epsilon-\rho\a_i =  r_i(\eps)$, $\epsilon-\sigma\a_j =
r_j(\eps)$, $\epsilon-(\rho+\sigma)\a_i-\sigma\a_j =r_ir_j(\eps)$,
$\epsilon-\rho\a_i-(\rho+\sigma)\a_j= r_jr_i(\eps)$,
$\epsilon-(\rho+\sigma)(\a_i+\a_j)=
r_jr_ir_j(\eps)=r_ir_jr_i(\eps)$. Note that $\rho = \sigma = 1$
would imply $(\a_j,\epsilon-\a_i) = 2$, whence
$\epsilon=\a_i+\a_j$.  This means $r_ir_jB=r_jB$.  To see this,
suppose $\a_i+\a_j\in B$.  Then $\a_i\in r_jB$ as it is
$r_j(\a_i+\a_j)$.  Now all other elements of $r_jB$ are orthogonal
to $\a_i$ as $\a_i$ is one of the elements.  Now $r_ir_jB=r_jB$.
So we can assume this does not occur.

By assumption, there are $\b,\c \in B$ such that $\b-\a_i\in r_i
B$ and $\c-\a_j\in r_j B$ and such that $\het(\b),\het(\c)$ are
minimal with respect to being moved by $\a_i$, $\a_j$
respectively. By symmetry, we may also assume that $\het(\b)\le
\het(\c)$.

If $\b=\c$, then $(\a_i,\b) = (\a_j,\c)=1$, a case that has been
excluded. Therefore, we may assume that $\b$ and $\c$ are
distinct. In particular, $(\a_j,\b)$ is $0$ or $-1$.

Suppose first that a $\b$ can be chosen so that $(\a_j,\b)=0$.
This is certainly the case if $\het(\b)<\het(\c)$. Now
$\b-\a_i-\a_j\in r_ir_jB$, so $\b-\a_i\in r_i B$ is a root of
smallest height moved by $r_j$ and so $r_jr_iB<r_i B$.  Recall
from our choice no root of height smaller than $\het(\b)$ is moved
by $r_i$.

Since $\b\in r_jB$ and $\b-\a_i \in r_ir_jB$, we have
$r_ir_jB<r_jB$ unless there is $\delta\in B$ with $\delta-\a_j\in
r_jB$, $\het(\delta-\a_j) = \het(\b)-1$, and $(\a_i,\delta-\a_j) =
-1$. But then the inner products show $\delta=-\a_i$ is not a
positive root.  Notice this shows $\delta-\a_j+a_i$ cannot be a
root.  Hence, indeed, $r_ir_jB<r_jB$.

Since $\b-\a_i-\a_j\in r_jr_ir_jB$ and $\b-\a_i\in r_ir_jB$, a
similar argument to the previous paragraph shows that $r_jr_ir_jB
< r_ir_jB$.

It remains to show $r_jr_ir_jB<r_jr_iB$. Both sides contain
$\b-\a_i-\a_j$ and $r_i$ does not lower or raise $r_jr_i\b$.  We
need to look at the $\delta$ in $B$ of height up to $\het(\c)$.
We know that for $\delta$ with $\het(\delta)<\het (\c)$ that
$(\delta,\a_j)=0$.  Looking at the equations above with $\sigma=0$
we see $r_i$ does not change
$r_jr_i\delta=\delta-\rho(\a_i+\a_j)$.  We also know that
$(\c,\a_j)=1$.  This means $\sigma =1$.  Using the equations again
with $\sigma=1$ and $\rho$ we must compare
$\c-\rho(\a_i+\a_j)-\a_j$ with $\c-\rho(\a_i+\a_j)-\a_j-\a_i $
which is lower.  In particular, $r_ir_jr_iB<r_jr_iB$.

Suppose then that $(\a_j,\b)=-1$.  This means in particular that
$\het(\b)=\het(\c)$. If $(\a_i,\c)=0$ we can use the argument
above.  We are left then with the case in which $(\a_i,\c)=1$,
$(\a_j,\b)=1$, $(\a_i,\b)=-1$, $(\a_j,\c)=-1$, and of course
$(\a_i,\a_j)=-1$.

\bigskip
This means $\c-\a_i$ and $\b-\a_j$ are positive roots of height
$\het(\b)-1$. But $(\c-\a_i , \b-\a_j) = 0+1+1-1 = 1$ so by
subtracting one root from the other, we should get another
positive root. As both roots are of the same height, this would
give a root of height $0$ which is not possible, proving this case
never arises.
\end{proof}

We showed during the proof that if $\a_i+\a_j\in B$ and $i\sim j$,
then $r_ir_jB=r_jB$.  This is case (iv) of Proposition
\ref{Prop:AdmOrd}.  The following lemma shows this is if and only
if.

\begin{lemma}\label{heighttworoot}
Suppose that $(\B,<)$ is a monoidal poset for $(W,R)$ for which
$r_ir_jB=r_jB$ with $i\sim j$.  If $r_iB<B$ and $r_jB<B$, then
$\a_i+\a_j\in B$.
\end{lemma}

\begin{proof}
Suppose $\a_i+\a_j$ is not in $B$.  Let $\b$ be an element of
smallest height moved in $B$ by $r_j$ for which $\b-\a_j$ is a
root.  Such a root exists because $r_jB<B$.  As $\a_i+\a_j$ is not
in $ B$, we know $\b-\a_j\neq \a_i$, and even $\a_i$ is not in
$r_jB$.  It follows as $r_ir_jB=r_jB$ that $\a_i\in (r_jB)^\perp$.
In particular $r_i(\b-\a_j)$ is in $r_jB$.  As $\b-\a_j\pm \a_i$
is not orthogonal to $\b-\a_j$ we must have $\b+\a_i$ a root.  Now
$r_j$ lowers $\b$ and $r_i$ raises $\b $.

As $r_iB<B$ there exists $\c$, an element of smallest height in
$B$ moved by $r_i$ for which $\c-\a_i$ is a root.  We know
$\het(\c)\leq \het(\b)$ as $r_i$ moves $\b$. Suppose
$(\c,\a_j)=0$.  Then $\c\in r_jB$ and $r_i\c=\c-\a_i$ is also in
$r_jB$.  This contradicts the hypothesis that elements of $r_iB$
are all orthogonal.  This implies $(\c,\a_j)=\pm 1$.  This in turn
means $\het(\b)\leq \het(\c)$ as $\het(\b)$ is the height of the
smallest element moved by $r_j$. Now we have $\het(\b)=\het(\c)$.
If $\c-\a_i$ and $\c-\a_j$ were both roots, an inner product
computation would show $(\c,\a_i+\a_j) =2$ so $\c=\a_i+\a_j$. This
means $\c+\a_j$ is a positive root, so in $r_jB$ we have $\b-\a_j$
and $\c+\a_j$ contradicting (i) of Proposition \ref{DnsystemIsAdm}.
\end{proof}

\bigskip
In order to address the monoid action later we will need some more
properties of this action in terms of lowering and raising.  We
begin with the case in which two different fundamental reflections
act the same on a member $B$ of $\B$.

Before we begin we need to examine the case in which some $B$ has
two indexes which raise it to the same $B'$.  In particular we
have

\begin{lemma}\label{goingupequal}
Suppose $B\in \B$ and $r_iB=r_kB>B$ with $k\ne i$.  If $\b $ is
the element of $B$ of smallest height moved by either $r_i$ or
$r_k$, then $\b+\a_i+\a_k$ is also in $B$.  Furthermore, $i\not
\sim k$.
\end{lemma}

\begin{proof}
Let $\b$ be an element of smallest height $B$ moved by either
$r_i$ or $r_k$. We know that all elements of smaller height are
not moved by $r_i$ and $r_k$.  Elements of the same height could
be moved by $r_i$ or $r_k$, but then the root would have to be
added.  Suppose $(\a_i,\b)=-1$, so $r_i\b  = \b+\a_i$. If
$(\a_k,\b)=0$,  then $\b\in r_kB=r_iB$ as is $\b+\a_i$ and so
$(\b,\b+\a_i)=2-1\ne 0$, which contradicts that elements of $r_iB$
are mutually orthogonal.  In particular $(\a_k,\b)=-1$ (for
otherwise, $(\a_k,\b) = 1$ and so $r_kB<B$).

If $i\sim k$, then $(\a_i,\a_k)=-1$, and so $(\a_k,\b+\a_i)=-2$,
which implies that $\b+\a_i=-\a_k$, contradicting that $\b+\a_i$
be a positive root.  This means $i\not \sim k$ which proves the
last part of the lemma.

Now by hypothesis $r_ir_kB=B$ and so $\b+\a_k+\a_i$ is in $B$
which proves the remainder of the lemma.
\end{proof}

Notice that if $\b$ and $\b+\a_i+\a_k$ are two roots in $B$ with
$(\a_i,\a_k)=0$, $(\b,\a_i)=(\b,\a_k)=-1$, the hypothesis of the
lemma is satisfied, and $r_i$ maps $\b$ to $\b+\a_i$ and
$\b+\a_i+\a_k$ to $\b+\a_k$.  Acting by $r_k$ has the same effect
except the order of the roots has been interchanged.

\begin{Lm}\label{Cor:AdmOrdProps}
Suppose $(\B,<)$ is a monoidal poset for $(W,R)$. Let $B\in\B$ and
let $i,j\in M$ and $\b,\c\in B$. Then the following assertions
hold.
\begin{enumerate}[(i)]
\item If $i\not\sim j$ and $r_ir_j B<r_i B<B$, then
$r_ir_jB<r_jB<B$.
\item If $i\not\sim j$, $B<r_i B$, $B<r_j B$,
and $r_iB\ne r_jB$, then $r_ir_j B>r_i B$ and $r_ir_jB>r_j B$.
\item If $i\sim j$, $B<r_i B$, and $B <r_j B$, then $r_i
B<r_jr_iB<r_ir_jr_iB$, and $r_j B<r_ir_j B<r_jr_ir_jB$. \item If
$i\sim j$ and $r_jr_ir_jB<r_jr_iB<r_i B <B$, then also
$r_jr_ir_jB<r_ir_j B<r_j B<B$.
\item If $\a_i\not\in B^\perp\cup B$, then either $r_iB<B$ or
$r_iB>B$.
\end{enumerate}
\end{Lm}

\begin{proof}
We can refer to Proposition \ref{Prop:AdmOrd} for the properties
of $(\B,<)$.

\nl(i).  If $r_j B = r_ir_j B$, then also $r_i B = B$, a
contradiction. Suppose $r_ir_j B>r_j B$. Then, by transitivity
$B>r_jB$. Also $B>r_i B$ by hypotheses.  Notice if $r_iB=r_jB$,
$r_ir_jB=r_i^2B=B$ but $r_ir_jB<B$.  Now $r_ir_j B <r_j B$ by
Proposition \ref{Prop:AdmOrd}(iii), a contradiction. Hence, by
Proposition \ref{Prop:AdmOrd}(i), $r_ir_j B<r_j B$.

If $r_j B = B$, then also $r_i B = r_ir_j B$, a contradiction.
Suppose $r_j B>B$. Then, by transitivity, $r_j B>r_ir_j B$.
Therefore, Proposition \ref{Prop:AdmOrd}(iii) gives $r_ir_j
B>r_jr_ir_jB = r_i B $, a contradiction. Hence, by Proposition
\ref{Prop:AdmOrd}(i), $r_j B<B$. But also $r_i B<B$, so
Proposition \ref{Prop:AdmOrd}(iii) gives $r_ir_j B<r_j B$ (and
$r_ir_jB<r_iB$).

\nl(ii).  If $r_ir_j B = r_i B$, then $r_j B = B$, a
contradiction.  If $r_ir_j B<r_i B$, then, by Proposition
\ref{Prop:AdmOrd}(iii) applied to $r_iB$ we have $r_j B < B$, a
contradiction. Hence by Proposition \ref{Prop:AdmOrd}(i), $r_ir_j
B>r_iB$. The proof of $r_ir_j B>r_j B$ is similar.

\nl(iii).  Suppose $r_ir_jB=r_jB$.  If $\a_i\in (r_jB)^\perp$,
then, as $r_j$ lowers $r_jB$, by Proposition \ref{Prop:AdmOrd}(ii)
$r_i$ lowers $r_jr_jB=B$ which is a contradiction.  This means
$\a_i\in r_jB$ and so $\a_i+\a_j\in B$.  Notice neither $\a_i$ nor
$\a_j$ are in $B$ as they are not orthogonal to $\a_i+\a_j$.  As
both $r_i$ and $r_j$ raise $B$, there must be $k$, $l$, with
$i\sim k$ and $j\sim l$ with $\a_k$ and $\a_l$ in $B$. Neither are
orthogonal to $\a_i+\a_j$ and this is impossible.  This means
$r_ir_jB\ne r_jB$.

Suppose $r_ir_jB<r_jB$.  We can't have $r_ir_jB=r_jr_jB=B$ by
Lemma \ref{goingupequal}.  Now Proposition \ref{Prop:AdmOrd}(iv)
gives $r_i B<B$, a contradiction. Hence $r_ir_j B>r_j B$.  The
roles of $i$ and $j$ are symmetric, so similarly we find
$r_ir_jB>r_i B$.

If $r_ir_jr_iB=r_ir_j B$ then $B = r_i B$, a contradiction.
Suppose $r_jr_ir_jB<r_ir_j B$. As also $r_j B<r_ir_j B$,
Proposition \ref{Prop:AdmOrd}(iv) gives $r_i B<B$, a
contradiction, because $\a_i+\a_j\in r_ir_j B$ would imply
$\a_j\in r_j B$ whence $r_jB=B$.

Similarly, it can be shown that $r_jr_ir_jB>r_ir_jB$.

\nl(iv). If $r_j B = B$, then $r_jr_ir_jB = r_jr_iB$, a
contradiction. If $r_j B<B$, then the result follows from
Proposition \ref{Prop:AdmOrd}(iv) because $\a_i+\a_j\in  B$ would
imply $\a_j\in r_i B$ whence $r_j r_iB=r_iB$.

Suppose therefore $r_j B>B$.  If $r_j B = r_ir_j B$, then
$r_jr_iB=r_jr_ir_jB$, a contradiction.  If $r_j B>r_ir_j B$, then
by Proposition \ref{Prop:AdmOrd}(iv) $r_jr_ir_jB>r_jr_iB$, a
contradiction because $\a_i+\a_j\in r_j B$ would imply $\a_i\in B$
whence $r_iB=B$.

Hence $r_j B<r_ir_j B$. But then by transitivity $r_ir_j
B>r_jr_ir_jB$, and, since $r_ir_j B>r_j B$, gives Proposition
\ref{Prop:AdmOrd}(iv) $r_jr_ir_jB>r_ir_jB$ (for otherwise
$\a_i+\a_j\in r_ir_j B$, implying $\a_j\in r_j B$ so $r_jB=B$), a
final contradiction.

\nl(v).  The hypotheses imply that there exists $\b\in B$ with
$(\a_i,\b)=\pm1$.  Then $r_i\b = \b\pm\a_i$, which is not
orthogonal to $\b$. As the elements of $B$ are orthogonal by
definition, $r_i\b$ does not belong to $B$, so $r_iB\ne B$, and
the conclusion follows from Proposition \ref{Prop:AdmOrd}(i).
\end{proof}

Pick $B_0$ a maximal element of $\B$.  This means $r_iB_0$ is
either $B_0$ or lowers $B_0$.  This is possible as $\B$ is finite.
We need more properties of the poset determined by $>$. To begin
with this we consider certain Weyl group elements, $w$, for which
$wB_0=B$ for a fixed element $B\in \B$. In particular we let
$w=r_{i_1}r_{i_2}\dots r_{i_s}$ be such that
$B_0>r_{i_s}B_0>r_{i_{s-1}}r_{i_s}B_0>\cdots >r_{i_2}r_{i_3}\cdots
r_{i_s}B_0>r_{i_1}r_{i_2}r_{i_3}\cdots r_{i_s}B_0=B$.  If there is
such an expression for $w$, then there is one of minimal length.
We let $\B'$ be the set of $B\in \B$ which are of this form.  We
will show that in fact $\B'=\B$.

\begin{Lm}\label{Ordergoingdown}  In the notation just above, $\B'=\B$.
\end{Lm}

\begin{proof}
Notice that $B_0$ is in $\B'$ by definition using $w$ the
identity.  Recall that $r_iB_0$ is either $B_0$ or lower.  In
particular nothing raises $B_0$. We show first that if $B\in \B'$
and $r_jB>B$ then $r_jB\in \B'$.  We prove this by induction on
the minimal length of a chain from $B_0$ to $wB$ which satisfies
the descending property of the definition of $\B'$.  In particular
$w=r_{i_1}r_{i_2}\dots i_s$ and
$B_0>r_{i_s}B_0>r_{i_{s-1}}r_{i_s}B_0>\cdots
>r_{i_2}r_{i_3}\cdots r_{i_s}B_0>r_{i_1}r_{i_2}r_{i_3}\cdots r_{i_s}B_0=B$.
We say this chain has length $s$, the length of $w$.  We in fact
show that there is a chain from $B_0$ to $r_jB$ of length less
than or equal to $s-1$.  We have seen that no $r_i$ raises $B_0$.
Suppose that $B=r_iB_0$ and $r_jB>B$.  If $r_jB=B_0$ the induction
assumption is true.  If $r_jB\neq B_0$ we can use Lemma
\ref{Cor:AdmOrdProps}(ii) or (iii) to see that $r_jr_iB>B_0$ a
contradiction.  In particular the induction assumption is true for
$s\le1$.

We can now assume $s\geq 2$.  Pick a $B$ with a chain of length
$s$ and assume the result is true for any $B'\in \B'$ with a
shorter chain length.  Suppose $r_jB$ is not in $\B'$ and
$r_jB>B$.  Notice $r_{i_1}B=r_{i_2}r_{i_3}\cdots r_{i_s}B_0>B$ by
the hypothesis.  Clearly $r_jB\neq r_{i_1}B$ as $r_{i_1}B $ is in
$\B'$ using the element $r_{i_2}r_{i_3}\cdots r_{i_s}$.  In
particular we can use Lemma \ref{Cor:AdmOrdProps}(ii) or (iii).
In either case $r_jr_{i_1}B>r_{i_1}B$ and by our choice of $s$ and
the induction assumption, $r_jr_{i_1}B$ is in $\B'$ and has a
chain of length at most $s-1$ from $B_0$ to it.

Suppose first $i_1\not \sim j$ and use Lemma
\ref{Cor:AdmOrdProps}(ii).  By the induction assumption there is a
chain down to $r_jr_{i_1}B$ of length at most $s-2$ and then by
multiplying by $r_{i_1}$ gives a chain down to $r_jB$ of length at
most $s-1$ and the induction gives $r_jB\in\B'$.

Suppose now $i_1\sim j$ and use Lemma \ref{Cor:AdmOrdProps}(iii).
Again $r_jr_{i_1}B$ is in $\B'$ by the induction hypothesis and
has a chain down to it of length at most $s-2$.  Using the
induction again, and the hypothesis of the minimality of $s$, we
see also $r_{i_1}r_jr_{i_i}B$ is in $\B'$ and has a chain to it of
length at most $s-3$.  Now using this as $r_jr_{i_1}r_jB$,
multiplying by $r_j$ and then by $r_{i_1}$ gives a chain to $r_jB$
of length at most $s-1$ and we are done with this part.

In particular, if $B\in \B'$ and $r_jB>B$, then $r_jB$ is in
$\B'$.  If $B\in \B'$ and $r_jB=B$ of course $r_jB\in \B'$.
Suppose $r_jB<B$.  Then the sequence to $B$ and then $r_jB$ gives
a sequence to $r_jB$ and $r_jB$ is in $\B'$.  We see that $\B'$ is
closed under the action of $W$ and as $\B$ is an orbit, $\B'=\B$.
\end{proof}

\begin{Cor}\label{uniquemax}
There is a unique maximal element $B_0$ in $\B$.
\end{Cor}

\begin{proof}  We have just shown that for every element $B$ in $\B$ except $B_0$ there is a sequence lowering to $B$ and
so $B_0$ is the only maximal element.
\end{proof}
\noindent See Example~\ref{A_n and D_n} for a listing of some of the $B_0$.
\medskip

This shows that each $B\in\B$ has a level associated with it,
namely the smallest $s$ for which $B$ can be obtained from $B_0$
as above with a Weyl group element $w$ of length $s$.  Namely the
smallest $s$ for which there is a reduced expression
$w=r_{i_1}r_{i_2}\cdots r_{i_s}$ with $wB_0=B$ for which
$B_0>r_{i_s}B_0>r_{i_{s-1}}r_{i_s}B_0>\cdots >r_{i_2}r_{i_3}\cdots
r_{i_s}B_0>B$.  In particlar $B_0$ has level $0$ and if
$r_jB_0<B_0$ it has level $1$.  The next lemma says that this $s$
is the shortest length of any word $w$ for which $wB_0=B$.

\begin{lemma}\label{shortestword}
Suppose $w$ is an element of $W$ of the smallest length for which
$wB_0=B$. Then this length, $s$, is the length of the shortest
word defining $B$ as an element of $\B'$.  In particular if the
word is $r_{i_1}r_{i_2}\cdots r_{i_s}$, then
$B_0>r_{i_s}B_0>r_{i_{s-1}}r_{i_s}\cdots>r_{i_1}r_{i_2}\cdots
r_{i_s}B_0=B$ and this is the shortest which does this.  It is
reduced.
\end{lemma}

\begin{proof}
Suppose $w$ is an element of $W$ for which $wB_0=B$ and for which
as in the definition of $\B'$, we have $w=r_{i_1}r_{i_2}\cdots
r_{i_s}$ and $r_{i_s}B_0>r_{i_{s-1}}r_{i_s}B_0>\cdots
>r_{i_1}r_{i_2}\cdots r_{i_s}B_0=B$ with this the shortest
possible.  Suppose $w'$ is any other Weyl group element with
$wB_0=B$. If $w'=r_{j_1}r_{j_2}\cdots r_{j_t}$ is a reduced
decomposition of length $t$, then $t$ is at most $s$ and we get a
sequence $B_0$, $r_{j_t}B_0$, $r_{j_{t-1}}r_{j_t}B_0$, $\cdots$,
$r_{j_1}r_{j_2}\cdots r_{j_t}B_0=B$.  If any of these differences
do not have the relation $>$ between them, the level of $B$ would
be strictly smaller than $s$, contradicting the minimality of $s$.
Hence, $t =s$ and the sequence corresponding to $w'$ is also a
chain. In particular, $w$ is reduced and any other reduced
expession gives a descending sequence of the same length.  This
shows there is a reduced word with this length taking $B_0$ to $B$
and any word doing this of shorter or the same length, has to be
descending at each step.  This proves the lemma.
\end{proof}

\begin{lemma}\label{lemma:dB}
Suppose that $(\B,<)$ is a monoidal poset for $(W,R)$.
\begin{enumerate}[(i)]
\item For each $B\in \B$ and each element $w\in W$ of minimal
length such that $B = wB_0$ and node $i$ of $M$ such that
$l(r_iw)<l(w)$, we have $r_iB>B$. \item For each $B\in \B$, if
$w,w'\in W$ are of minimal length such that $B = wB_0 = w'B_0$,
then $l(w) = l(w')$ and, for each node $i$ such that $r_iB>B$,
there is $w''\in W$ of length $l(w)$ such that $B= w''B_0$ and
$l(r_iw'') < l(w'')$.
\end{enumerate}
\end{lemma}

\begin{proof}
For (i) we use the characterization in Lemma \ref{shortestword}
and realize that any of the equivalent expressions also give a
descending sequence.  In particular if $l(r_iw)<l(w)$, an
equivalent word can be chosen to start with $r_i$ and so $r_iB$ is
one step above $B$ in the chain to $B$ from this word and so
$r_iB>B$.

For (ii) again use Lemma \ref{shortestword} and so $l(w)=l(w')$.
If $r_iB>B$ for some $i$, there is a sequence from $r_iB$ to
$B_0$.  If $w'B_0=r_iB$ accomplishes this in the minimal number of
steps, $w''=s_iw'$ satisfies the conclusion of the lemma.
\end{proof}

\section{The positive monoid}\label{def:h}
We now turn our attention to the Artin group $A$ associated with
the Coxeter system $(W,R)$. We recall that the defining
presentation of $A$ has generators $s_i$ corresponding to the
fundamental reflections $r_i\in R$ and {\em braid relations}
$s_is_js_i = s_js_is_j$ if $i\sim j$ and $s_is_j = s_js_i$ if
$i\not\sim j$.  The monoid $A^+$ given by the same presentation is
known (\cite{paris}) to embed in $A$.  For each admissible
$W$-orbit of a set of mutually commuting reflections, we shall
construct a linear representation of $A^+$. To this end, we need a
special element $h_{B,i}$ of $A^+$ for each pair $(B,i)$
consisting of a set $B$ of mutually commuting reflections and a
node $i$ of $M$ whose reflection $r_i$ does not belong to $B$ but
commutes with each element of $B$. As in the previous section, we
shall represent reflections by positive roots.

We now define the elements $h_{B,i}$.
As in \cite{CGW} we do this by defining reduced words $v_{B,i}\in
A$ and letting $h_{B,i}={v_{B,i}^{-1}}s_i{v_{B,i}}$. Later we
shall consider the image of these elements in a certain Hecke
algebra.

We make definitions of $v_{B,i}$ which depend on certain chains
from $B_0$ to $B$ and show in an early lemma that conjugating
$s_i$ by any of them gives the same element.  Furthermore, this
element corresponds to a fundamental generator of $A$ commuting
with every reflection having its positive root in $B_0$.

\begin{Def}\label{defw} Suppose $(\B,<)$ is a monoidal poset of $(W,R)$,
with maximal element $B_0$. Let $(B,i)$ be a pair with $B\in\B$
and $i$ a node of $M$ such that $\a_i\in B^\perp$.

Choose a node $j$ of $M$ with $r_jB>B$.  If $j\not \sim i$ let
$v_{B,i}=s_jv_{r_jB,i}$ and if $i\sim j$ let
$v_{B,i}=s_js_iv_{r_ir_jB,j}$. We define $v_{B_0,i}$ as the
identity.

Furthermore, set $h_{B,i}= {v_{B,i}}^{-1} s_i {v_{B,i}}$.
\end{Def}

Notice this definition makes sense as a nondeterministic algorithm
assigning an element of $A$ to each pair $(B,i)$ as specified
because
\begin{itemize}
\item if $\a_i\in B^\perp$ and $i\not \sim j$, then $\a_i\in\
(r_jB)^\perp$; \item  If $i\sim j$, then $\a_i+\a_j \in
(r_jB)^\perp$ and $\a_j\in (r_ir_jB)^\perp$.
\end{itemize}
By Lemma \ref{shortestword}, $v_{B,i}$ will be a reduced
expression whose length is the length of a chain from $B$ to
$B_0$. The elements $v_{B,i}$ are not uniquely determined, but we
will show that the elements $h_{B,i}$ are.

\medskip

\begin{lemma}\label{mainlemma}
For $\B$ and $(B,i)$ as in Definition \ref{defw}, suppose that
$v_{B,i}$ and $v_{B,i}'$ both satisfy Definition \ref{defw}.  Then
the elements $h_{B,i}$ of $A$ defined by each are the same, i.e.,
$h_{B,i}={v_{B,i}}^{-1} s_i {v_{B,i}}= {v_{B,i}}'^{-1} s_i
{v_{B,i}}'$.

Furthermore, each $h_{B,i}$ is a fundamental generator $s_j$ of
$A$ whose root $\a_j$ is orthogonal to every root of $B_0$.
\end{lemma}

\begin{proof}
We use induction on the height from $B_0$.  The case of $B_0$ is
trivial.

We first dispense with the case in which $r_jB=r_{j'}B$.  We know
from Lemma \ref{goingupequal} that in $B$ there are two elements
$\b$ and $\b+\a_j+\a_{j'}$ with $j\not \sim j'$. As $\a_i\in
B^\perp$ we know $(\b,\a_i)=0$ and also
$(\b+\a_j+\a_{j'},\a_i)=0$.  It follows that
$(\a_j,\a_i)=(\a_{j'},\a_i)=0$ as the inner products of
fundamental roots are $0$ or $-1$.  In particular using $r_j$ we
get $v_{B,i}=s_jv_{r_jB,i}$, and
${v_{B,i}}^{-1}s_i{v_{B,i}}={v_{r_jB,i}}^{-1}s_j^{-1}s_is_j{v_{B,i}}$.
As $s_j^{-1}s_is_j=s_i$ this is $
{v_{r_jB,i}}^{-1}s_i{v_{B,i}}=h_{r_jB,i}$.  The same is true for
$r_{j'}$ and we are assuming $r_jB=r_{j'}B$.  Now we can use
induction.

We next suppose $r_jB>B$ and $r_{j'}B>B$ with $r_jB\ne r_{j'}B$.
There will be two cases depending on whether $j\not \sim j'$ or
$j\sim j'$.  Suppose first $j\not \sim j'$.  We use Lemma
\ref{Cor:AdmOrdProps}(ii) to see that $r_jr_{j'}B>r_{j'}B$ and
$r_{j'}r_jB>r_jB.$  Suppose first $i\not \sim j$ and $i\not \sim
j'$.  For the chain starting with $r_j$ we can follow it with
$r_{j'}$ and if we start with $r_{j'}$ we can follow it with
$r_j$.  In each case with these choices we get
$v_{B,i}=s_js_{j'}v_{r_jr_{j'}B,i}$ as $s_js_{j'}=s_{j'}s_j$ and
$\a_i\in (r_jB)^\perp$ and $\a_i\in r_{j'}B)^\perp$.  The
induction is used for $v_{r_jB,i}$ and for $v_{r_{j'}B,i}$ in
order to take the chain we have chosen and then also for
$v_{r_jr_{j'}B,i}$.  In each case we get $h_{r_jr_{j'}B,i}$.

Suppose next that $i \sim j$ but $i\not \sim j'$.  Using the chain
for $r_j$ we get $B<r_jB<r_ir_jB$ by Proposition
\ref{Prop:AdmOrd}(ii).  As above by Lemma
\ref{Cor:AdmOrdProps}(ii) we get $r_jr_{j'}B>r_{j'}B$ and now
again by Proposition \ref{Prop:AdmOrd}(ii) using $\a_i\in
(r_{j'}B)^\perp$ we get $r_ir_jr_{j'}B>r_jr_{j'}B$.  Now for the
$r_{j'}$ chain continue through $r_j$ and then $r_i$ to reach
$v_{B,i}=s_{j'}s_js_iv_{r_ir_jr_{j'}B,j}$.  Through the $r_j$
chain which goes through $r_ir_jB$ add $r_{j'}$ for which $j'\not
\sim i$.  Here we get $v_{B,i}=s_js_is_{j'}v_{r_{j'}r_ir_jB,j}$.
Again use induction at all the levels to get the needed result.
Notice $r_ir_jB\ne r_{j'}r_jB$ as $r_ir_jr_{j'}B>r_jr_{j'}B$ as
above and so $r_i $ raises $r_jr_{j'}B$.

The final case in which $j\not \sim j'$ is with $j\sim i\sim j'$,
see Figure 1.  For this we again use Lemma
\ref{Cor:AdmOrdProps}(ii) and (iii) and Proposition
\ref{Prop:AdmOrd}.  In particular $r_jB > B$ and $r_ir_jB > r_jB$.
Also as $r_jB\ne r_{j'}B$ we have $r_{j'}r_jB>r_jB$.  Now by Lemma
\ref{Cor:AdmOrdProps}(iii) we have $r_ir_{j'}r_ir_jB
>r_{j'}r_ir_jB >r_ir_jB$. We know $\a_j\in (r_ir_jB)^\perp $ and
also in $(r_{j'}r_ir_jB)^\perp $ as $j\not \sim j'$.  Now using
Proposition \ref{Prop:AdmOrd}(ii) we see
$r_jr_ir_{j'}r_ir_jB>r_ir_{j'}r_ir_jB$.  Notice $r_ir_jB\ne
r_{j'}r_jB$ by Lemma \ref{goingupequal} as $i\sim j'$.  Following
the trail of $\a_i $ we see it is in
$(r_jr_ir_{j'}r_ir_jB)^\perp$.  Following this chain after
$r_ir_j$ and using induction we see
$v_{B,i}=s_js_is_{j'}s_is_jv_{r_jr_ir_{j'}r_ir_jB,i}$.  Going up
through $r_ir_{j'}$ gives the same result as
$s_js_is_{j'}s_is_j=s_js_{j'}s_is_{j'}s_j$ is similar to
$s_{j'}s_is_js_is_{j'}=s_{j'}s_js_is_js_{j'}$.  In particular the
result is true again using induction at all the higher levels.

\unitlength .7mm
\begin{picture}(100,120)(-35,05)
\linethickness{0.3mm} \put(60,93.54){\circle{7.07}}
\linethickness{0.3mm} \put(30,93.54){\circle{7.07}}
\linethickness{0.3mm} \put(45,113.54){\circle{7.07}}
\linethickness{0.3mm}
\multiput(47.5,111.04)(0.12,-0.18){83}{\line(0,-1){0.18}}
\linethickness{0.3mm}
\multiput(32.5,96.04)(0.12,0.18){83}{\line(0,1){0.18}}
\put(15,73.54){\circle{7.07}}

\linethickness{0.3mm}
\multiput(17.5,76.04)(0.12,0.18){83}{\line(0,1){0.18}}
\linethickness{0.3mm}
\multiput(47.5,76.04)(0.12,0.18){83}{\line(0,1){0.18}}
\linethickness{0.3mm}
\multiput(62.5,91.04)(0.12,-0.18){83}{\line(0,-1){0.18}}
\linethickness{0.3mm} \put(75,73.54){\circle{7.07}}

\linethickness{0.3mm}
\multiput(32.5,91.04)(0.12,-0.18){83}{\line(0,-1){0.18}}
\linethickness{0.3mm} \put(45,73.54){\circle{7.07}}

\linethickness{0.3mm} \put(75,50){\line(0,1){20}}
\linethickness{0.3mm} \put(45,50){\line(0,1){20}}
\linethickness{0.3mm} \put(15,50){\line(0,1){20}}
\linethickness{0.3mm} \put(30,26.46){\circle{7.07}}

\linethickness{0.3mm} \put(60,26.46){\circle{7.07}}

\linethickness{0.3mm} \put(45,6.46){\circle{7.07}}

\linethickness{0.3mm}
\multiput(32.5,23.96)(0.12,-0.18){83}{\line(0,-1){0.18}}
\linethickness{0.3mm}
\multiput(47.5,8.96)(0.12,0.18){83}{\line(0,1){0.18}}

\put(75,46.46){\circle{7.07}} \linethickness{0.3mm}
\multiput(62.5,28.96)(0.12,0.18){83}{\line(0,1){0.18}}
\linethickness{0.3mm}
\multiput(32.5,28.96)(0.12,0.18){83}{\line(0,1){0.18}}
\linethickness{0.3mm}
\multiput(17.5,43.96)(0.12,-0.18){83}{\line(0,-1){0.18}}
\linethickness{0.3mm} \put(15,46.46){\circle{7.07}}

\linethickness{0.3mm}
\multiput(47.5,43.96)(0.12,-0.18){83}{\line(0,-1){0.18}}
\linethickness{0.3mm} \put(45,46.46){\circle{7.07}}

\put(45,0){\makebox(0,0)[cc]{$B$}}

\put(78,60){\makebox(0,0)[cc]{$j$}}

\put(48,60){\makebox(0,0)[cc]{$i$}}

\put(18,60){\makebox(0,0)[cc]{${j'}$}}

\put(42,103){\makebox(0,0)[cc]{$j$}}

\put(50,103){\makebox(0,0)[cc]{${j'}$}}

\put(50,17){\makebox(0,0)[cc]{${j'}$}}

\put(42,17){\makebox(0,0)[cc]{$j$}}

\put(20,35){\makebox(0,0)[cc]{$i$}}

\put(42,35){\makebox(0,0)[cc]{${j'}$}}

\put(50,35){\makebox(0,0)[cc]{$j$}}

\put(70,35){\makebox(0,0)[cc]{$i$}}

\put(70,85){\makebox(0,0)[cc]{$i$}}

\put(50,85){\makebox(0,0)[cc]{$j$}}

\put(42,85){\makebox(0,0)[cc]{${j'}$}}

\put(20,85){\makebox(0,0)[cc]{$i$}}
\end{picture}
\begin{center}Figure 1
\end{center}

\medskip
We now consider the cases in which $j\sim j'$.  As before, the
easiest case is when $i\not \sim j$ and $i \not \sim j'$.  In this
case as before use Lemma \ref{Cor:AdmOrdProps}(iii) to obtain
$r_{j'}r_jr_{j'}B>r_jr_{j'}B>r_{j'}B$.  As $\a_i$ is orthogonal to
$\a_j$ and $\a_{j'}$ we obtain
$v_{B,i}=s_{j'}s_js_{j'}v_{r_{j'}r_jr_{j'}B,i}$.  The same result
holds for the other order and the result follows by induction.

The only other possibility is $i\sim j$ and $i\not \sim j'$ as $i$
could not be adjacent to both $j$ and $j'$ (for otherwise there
would be a triangle in the Dynkin diagram).  For this we use the
familiar six sides diagram generated by $r_j$ and $r_{j'}$ using
Lemma \ref{Cor:AdmOrdProps}(iii), see Figure 2. At $r_jB$ we may
also act by $r_i$ which we know raises $r_jB$.  It is clear that
$r_{j'}r_jB\ne r_ir_jB$ as $\a_j$ is in $(r_ir_jB)^\perp $ and so
$r_j$ could not raise (or even lower) it.  We can proceed by Lemma
\ref{Cor:AdmOrdProps}(ii) to $r_{j'}r_jr_ir_jB$.  Notice
$r_jr_{j'}r_jB \ne r_ir_{j'}r_jB$ by Lemma \ref{goingupequal} as
$i\sim j$.  Now proceed by Lemma \ref{Cor:AdmOrdProps}(iii) to
$r_ir_jr_ir_{j'}r_jB=r_ir_jr_{j'}r_ir_jB$.  By following the
perpendicularities we see $\a_{j'}\in
(r_ir_jr_{j'}r_ir_jB)^\perp$.  Using this chain which starts with
$r_j$ and continues with $r_i$, we find
$v_{B,i}=s_js_is_{j'}s_js_iv_{s_ir_jr_{j'}r_ir_jB,j'}$.  Using the
other direction starting with $r_{j'}$ then $r_j$, then $r_i$, we
can continue with $r_{j'}$ and $r_j$ to get
$r_jr_{j'}r_ir_jr_{j'}B$ and conclude using this direction
$v_{B,i}=s_{j'}s_js_is_{j'}s_jv_{r_jr_{j'}r_ir_jr_{j'}B,j'}$.  At
the juncture $r_jr_{j'}B$ we act by $r_i$ or by $r_{j'}$.  These
two could not be equal as again $\a_j$ is in
$(r_ir_jr_{j'}B)^\perp$ and so $r_j$ could not move it.  However,
if they were equal, it lowers it to $r_{j'}r_jB$.  These words are
equivalent and we can use induction as usual for the last time.

\unitlength .7mm
\begin{picture}(170,107)(-10,0)
\put(45,10){\makebox(0,0)[cc]{}}

\put(50,0){\makebox(0,0)[cc]{$B$}}

\put(80,20){\makebox(0,0)[cc]{}}

\put(99,79){\makebox(0,0)[cc]{$j$}}

\put(63,95){\makebox(0,0)[cc]{$j$}}

\put(59,54){\makebox(0,0)[cc]{$j$}}

\put(60,18){\makebox(0,0)[cc]{$j$}}

\put(31,37){\makebox(0,0)[cc]{$j$}}

\put(87,55){\makebox(0,0)[cc]{$i$}}

\put(97,35){\makebox(0,0)[cc]{$i$}}

\put(53,78){\makebox(0,0)[cc]{$i$}}

\put(22,60){\makebox(0,0)[cc]{$i$}}

\put(89,98){\makebox(0,0)[cc]{$i$}}

\put(108,60){\makebox(0,0)[cc]{${j'}$}}

\put(68,36){\makebox(0,0)[cc]{${j'}$}}

\put(41,17){\makebox(0,0)[cc]{${j'}$}}

\put(27,79){\makebox(0,0)[cc]{${j'}$}}

\put(41,54){\makebox(0,0)[cc]{${j'}$}}

\linethickness{0.3mm}
\multiput(53,61)(0.24,-0.12){75}{\line(1,0){0.24}}
\linethickness{0.3mm}
\multiput(29,21)(0.18,-0.12){99}{\line(1,0){0.18}}
\linethickness{0.3mm} \put(73.53,23.47){\circle{7.06}}

\linethickness{0.3mm} \put(73,27){\line(0,1){19}}
\linethickness{0.3mm}
\multiput(53.42,8.71)(0.17,0.12){102}{\line(1,0){0.17}}
\linethickness{0.3mm} \put(50,7.54){\circle{7.07}}

\linethickness{0.3mm} \put(73.53,49.53){\circle{7.06}}

\linethickness{0.3mm} \put(26.47,49.53){\circle{7.06}}

\linethickness{0.3mm}
\multiput(29,52)(0.25,0.12){75}{\line(1,0){0.25}}
\linethickness{0.3mm} \put(50.46,63.54){\circle{7.07}}

\linethickness{0.3mm} \put(26.53,23.53){\circle{7.06}}

\linethickness{0.3mm} \put(27,27){\line(0,1){19}}
\linethickness{0.3mm} \put(113.54,49.46){\circle{7.07}}

\linethickness{0.3mm}
\multiput(76,26)(0.2,0.12){175}{\line(1,0){0.2}}
\linethickness{0.3mm}
\multiput(100,62)(0.13,-0.12){83}{\line(1,0){0.13}}
\linethickness{0.3mm} \put(13.53,63.53){\circle{7.06}}

\linethickness{0.3mm}
\multiput(16,66)(0.18,0.12){175}{\line(1,0){0.18}}
\linethickness{0.3mm}
\multiput(16,61)(0.12,-0.13){67}{\line(0,-1){0.13}}
\linethickness{0.3mm}
\multiput(76.06,101.29)(0.24,-0.12){75}{\line(1,0){0.24}}
\linethickness{0.3mm} \put(96.6,63.76){\circle{7.06}}

\linethickness{0.3mm} \put(96.06,67.29){\line(0,1){19}}
\linethickness{0.3mm}
\multiput(76,52)(0.23,0.12){77}{\line(1,0){0.23}}
\linethickness{0.3mm} \put(96.6,89.82){\circle{7.06}}
\linethickness{0.3mm} \put(49.53,89.82){\circle{7.06}}
\linethickness{0.3mm}
\multiput(52.06,92.29)(0.25,0.12){75}{\line(1,0){0.25}}
\linethickness{0.3mm} \put(73.53,103.83){\circle{7.07}}
\linethickness{0.3mm} \put(50.06,67.29){\line(0,1){19}}
\end{picture}
\begin{center}{Figure 2}
\end{center}

\medskip

This finishes all cases and shows the words have the same effect
under conjugation on $s_i$.
\end{proof}

We finish this section by exhibiting relations that hold for the
$h_{B,i}$. Since we are actually interested in their images in the
Hecke algebra $H$ of type $M$ under the natural morphism
$\Q(m)[A]\to H$, we phrase the result in terms of elements of this
algebra.

\begin{Prop}\label{prop:hBi}
Suppose that $(\B,<)$ is a monoidal poset with maximal element
$B_0$. Let $C$ be the set of nodes of $M$ such that $\a_i$ is
orthogonal to $B_0$ and denote $Z$ the Hecke algebra over $\Q(m)$
of the type $C$.  Then the images of the elements $h_{B,i}\in A$
in the Hecke algebra of type $M$ under the natural projection from
the group algebra of $A$ over $\Q(m)$ actually are fundamental
generators of $Z$ and satisfy the following properties.
\begin{enumerate}[(i)]
\item $h_{B,i}^2 = 1- mh_{B,i}$. \item $h_{B,i}h_{B,j} =
h_{B,j}h_{B,i}$ if $i\not\sim j$. \item $h_{B,i}h_{B,j}h_{B,i} =
h_{B,j}h_{B,i}h_{B,j}$ if $i\sim j$. \item $h_{r_jB,i} = h_{B,i}$
if $i\not\sim j$. \item $h_{r_ir_j B,j} = h_{B,i}$ if $i\sim j$
and $(\a_j,B)\ne 0$.
\end{enumerate}
\end{Prop}

\begin{proof}
By \cite{crisp} we can identify the Hecke algebra of type $C$ with
the subalgebra of the Hecke algebra $H$ generated by the $s_i$ for
$i\in C$.  As described above we define
$h_{B,i}={v_{B,i}}^{-1}s_i{v_{B,i}}$ where we consider this
element in the Hecke algebra. By Lemma \ref{mainlemma}, it is a
fundamental generator of $Z$.

\nl(i). This clearly follows from the quadratic Hecke algebra
relations we are assuming.

\nl(ii). Assume first $r_kB>B$ and both $\a_i$ and $\a_j$ are
orthogonal to $\a_k$.  We are assuming here $i\not \sim j$.  Then
we can take $v_{B,i}=s_kv_{r_kB,i}$ and $v_{B,j}=s_kv_{r_kB,j}$.
Now $h_{B,i}=({v_{r_kB,i}})^{-1}s_k^{-1}s_is_k{v_{r_kB,i}}$.  This
is $h_{r_kB,i}$ and we can use induction.

Suppose $i\sim k$ but $j\not \sim k$.  Then we can take
$v_{B,i}=s_ks_iv_{r_ir_kB,k}$ and we can take
$v_{B,j}=s_ks_iv_{r_ir_kB,j}$. Then $h_{B,j}=h_{r_ir_kB,j}$ and
$h_{B,i}=h_{r_ir_kB,k}$.  Now as above we can again use induction.

The final case with $i\not \sim j$ is when $i\sim k\sim j$.  Now
$v_{B,i}=s_ks_iv_{r_ir_kB,k}$ and $v_{B,j}=s_ks_jv_{r_jr_kB,k}$.
Suppose that $r_ir_kB=r_jr_kB$.  If so $v_{B,i}=s_ks_iw'$ with
$v_{B,j}=s_ks_jw'$ and $w'=v_{r_ir_kB,k}$.  Now
${v_{B,i}}^{-1}s_i{v_{B,i}}={w'^{-1}}s_i^{-1}s_k^{-1}s_is_ks_i{w'}
={w'^{-1}}s_k{w'}$. Doing the same with $r_j$ gives the same thing
and so they commute.  This means $r_jr_kB\neq r_ir_kB$ and we can
use Lemma \ref{Cor:AdmOrdProps}(ii) to get $r_ir_jr_kB>r_ir_kB$
and $r_ir_jr_kB>r_jr_kB$.  Now applying this with $v_{B,i}$ gives
$v_{B,i}=s_ks_js_is_kv_{r_kr_ir_jr_kB,j}$ and
$v_{B,j}=s_ks_is_js_kv_{r_kr_jr_ir_kB,i}$.  Let
$B'=r_kr_ir_jr_kB$.  Now $h_{B,i}=h_{B',j}$ and
$h_{B,j}=h_{B',i}$.  Now use induction.

\nl(iii). Suppose $i\sim j$.  We wish to show
$h_{B,i}h_{B,j}h_{B,i}=h_{B,j}h_{B,i}h_{B,j}$.  Suppose first
$k\not \sim i$ and $k\not \sim j$.  In this case
$v_{B,i}=s_kv_{r_kB,i}$ and $v_{B,j}=s_kv_{r_kB,j}$.  This means
$h_{B,i}=h_{r_kB,i}$ and $h_{B,j}=h_{r_kB,j}$.  Now use induction.

We are left with the case where $k\sim i\sim j$. Then $j\not \sim
k$ as there are no triangles in the Dynkin diagram.  Notice on the
chain from $\a_i$ we start with $r_k$, apply $r_i$ and can then if
we wish add $r_j$ provided $r_j$ raises $r_ir_kB$.  The chain from
$\a_j$ is $r_k$ which fixes $\a_j$, and then we can continue with
$r_i$ and then $r_j$ which forces $r_jr_ir_kB>r_ir_kB$ by
Proposition \ref{Prop:AdmOrd}(ii).  Now
$v_{B,i}=s_ks_is_jv_{r_jr_ir_kB,k}$ and
$v_{B,j}=s_ks_is_jv_{r_jr_ir_kB,i}$.  Now check that if
$B'=r_jr_ir_kB$ that $h_{B,i}=h_{B',k}$ and $h_{B,j}=h_{B',i}$.
Now use induction.

\nl(iv). Suppose $r_jB>B$.  Then $v_{B,i}=s_jv_{r_jB,i}$.  Now
conjugating $r_i$ by ${v_{B,i}}$ has the same effect as
conjugating ${v_{r_jB,i}}$ as $s_j^{-1}s_is_j=s_i$.  If $r_jB<B$,
use the same argument on $r_jB$ which is raised by $r_j$.

\nl(v). Assume first that $r_jB>B$.  Then
$v_{B,i}=s_js_iv_{r_ir_jB,j}$. Notice
$s_i^{-1}s_j^{-1}s_is_js_i=s_j$ and conjugating $s_i$ by
${v_{B,i}}$ has the same effect as conjugating $s_j$ by
${v_{r_ir_jB,j}}$ and the result follows.  If $r_jB<B$, then
$r_ir_jB<jB$ by Proposition \ref{Prop:AdmOrd}(ii). Now apply the
above to $r_ir_jB$.  As $(\a_j,B)\ne 0$, we know $r_jB\ne B$ by
Lemma \ref{Cor:AdmOrdProps}(v).

All cases have been completed.
\end{proof}

\begin{Remark}\label{rmk:w-not-all-paths}
\rm For the definition of $v_{B,i}$ we have used chains (and their
labels) from $B$ to $B_0$ depending on $i$.  In particular for
$r_jB>B $ and $i\not \sim j$ we use $s_jv_{r_jB,i}$ and for $j\sim
i$ we use $s_js_iv_{r_ir_jB,i}$.  If we were to use just any chain
we would not get this unique element without some further work.
For instance, if $M = {\rm D}_5$ and $B=\{
\eps_3+\eps_4,\eps_1+\eps_2\}$, both $\a_1$ and $\a_3$ are in
$B^\perp$ and $r_2r_1$ and $r_2r_3$ both take $B$ to $B_0=\{
\eps_1+\eps_4,\eps_2+\eps_3 \}$.  If we use the definition here,
with $v_{B,1}=s_2s_1$, we find $h_{B,1}=(s_2s_1)^{-1}s_1(s_2s_1) =
s_2$.  However, if we would use $v_{B,1}=s_2s_3$, corresponding to
a non-admitted chain, we find $s_3^{-1}s_2^{-1}s_1s_2s_3$ instead
of $s_2$ and we would need a proper quotient of the Hecke algebra
for $h_{B,1}$ to be well defined.
\end{Remark}

\begin{Cor}\label{cor:hB}
Let $(\B,<)$ be a monoidal poset. Retain the notation of the previous
proposition. Denote $C$ the set of all nodes $j$ of $M$ such that
$(\a_j,B_0) = 0$ and $Z$ the Hecke algebra whose type is the
diagram $M$ restricted to $C$. Then, for each node $j$ in $C$,
there is a minimal element $B$ of $(\B,<)$ and a node $k$ of $M$
such that $(\a_k,B) = 0$ and $h_{B,k} = s_j$, the image of the
fundamental generator of $A$ in $Z$.
\end{Cor}

\begin{proof}
The following proof is similar to the one of Lemma 3.8 of \cite{CGW}.
Let $j$ be a node of $C$. Then $h_{B_0,j} = s_j$. Let $B\in \B$ be
minimal such that there exists a node $k$ with $(\a_k,B) = 0$ and
$h_{B,k} = s_j$.  Suppose there is a node $i$ such that $r_i B<B$.
If $i\not\sim k$ then by Proposition \ref{prop:hBi}(iv) $h_{r_i
B,k} = h_{B,k} = s_j$.  If $i\sim k$ then by Proposition
\ref{prop:hBi}(v) $h_{r_kr_i B,i}=h_{B,k}=s_j$ and by Proposition
\ref{Prop:AdmOrd}(ii), $r_kr_i B<B$.  Both cases contradict the
minimal choice of $B$, so $B$ must be a minimal element of
$(\B,<)$.
\end{proof}

\begin{Example} \label{A_n and D_n}
Suppose $M$ is a connected simply laced diagram. Then the type of
$C$ as defined in Corollary \ref{cor:hB} is given in
Table~\ref{tab:type}. We deal with two series in particular.

If $M=\A_{n-1} $ and  $\B$ is the $W$-orbit
of $\{\a_1,\a_3,\ldots,\a_{2p-1}\}$,  then
\begin{eqnarray*}
B_0 & =
&\{\eps_1-\eps_{n-p+1},\eps_2-\eps_{n-p+2},\ldots,\eps_p-\eps_{n}\}
\mbox{ and}\\
C &=& \{\alp_{p+1},\alp_{p+2},\ldots,\alp_{n-p-1}\}.
\end{eqnarray*}
Therefore, the Hecke algebra $Z$ is of type $\A_{n-2p-1}$.

If $M=\D_n$ and $\B $ is the $W$-orbit of
$\{\a_1,\a_3,\ldots,\a_{2p-1}\}$,  then
\begin{eqnarray*}
B_0 & =
&\{\eps_1+\eps_{2p},\eps_2+\eps_{2p-1},\ldots,\eps_p+\eps_{p+1}\}
\mbox{ and }\\
C &=& \{\alp_p,\alp_{2p+1},\alp_{2p+2},\ldots,\alp_n\}.
\end{eqnarray*}
The Hecke algebra $Z$ has type $\A_1\D_{n-2p}$ (where $\D_1$ is
empty and $\D_2=\A_1\A_1$).
\end{Example}

\section{The Monoid Action}
Let $\B$ be an admissible $W$-orbit of sets of mutually orthogonal
positive roots, let $(\B,<)$ be the corresponding monoidal poset
(cf.\ Proposition \ref{Prop:AdmOrd}), let $B_0$ be the maximal
element of $(B,<)$ (cf.\ Corollary \ref{uniquemax}), and let $C$
be the set of nodes $i$ of $M$ with $\a_i\in B_0^\perp$. As before
(Proposition \ref{prop:hBi}), $Z$ is the Hecke algebra over
$\Q(m)$ of type $C$. These are listed in Table~$1$ under column $C$.
In analogy to the developments in \cite{CGW}
we define a free right $Z$-module $V$ with basis $x_B$ indexed by
the elements $B$ of $\B$.  By Lemma \ref{mainlemma} the linear
transformations $\tau_i$ of (\ref{taui}) are completely
determined.  We are ready to prove the main theorem.

\begin{proof}[Proof of Theorem \ref{th:main}]
Let $M$ be connected (see a remark following the theorem). We need
to show that the braid relations hold for $\tau_i$ and $\tau_j$, that
is, they commute if $i\not \sim j$ and
$\tau_i\tau_j\tau_i=\tau_j\tau_i\tau_j$ if $i\sim j$.

Take $B\in \B$. By linearity, it suffices to check the actions on
$x_B$.  We first dispense with the case in which either
$\tau_ix_B$ or $\tau_jx_B$ is $0$.  This happens if $B$ contains
$\a_i$ or $\a_j$.  If both roots are in $B$ both images are $0$
and the relations hold.

Suppose then that $\a_i$ is in $B$ but $\a_j$ is not in $B$.
Consider first the case in which $i \not \sim j$.  Then
$\tau_ix_B=0$ and so $\tau_j\tau_ix_B=0$.  Now $\tau_jx_B$ is in
the span of $x_B$ and $x_{r_jB}$. Notice as $(\a_i,\a_j)=0$ that
$\a_i$ is in $r_jB$ as well as $B$ and so $\tau_i\tau_jx_B=0$
also.  Suppose $i\sim j$.  Clearly $\t_i\t_j\t_ix_B=0$ as
$\t_ix_B=0$.  As $\a_i\in B$, the root $\a_i+\a_j$ belongs to
$r_jB$.  Also $r_j$ raises $B$ as a height one element, $\a_i$,
becomes height $2$.  This means $\t_jx_B=x_{r_iB}-mx_B$.  If $r_i$
lowers $r_jB$, $\t_ix_{r_iB}=x_{r_ir_jB}$. But $r_ir_jB$ contains
$r_i(\a_i+\a_j)=\a_j$, so $\t_jx_{r_ir_jB}=0$.  Also $\t_ix_B=0$
as $\a_i\in B$.  This proves the result unless $\t_i$ raises
$\t_jB$.  We know $\t_i$ takes the root $\a_i+\a_j$ to $\a_j$ and
so lowers a root of height $2$.  The only way $r_i$ could raise
$r_jB$ is if $r_jB$ contained an $\a_k$ with $k\sim i$.  This
would be $r_j\b$ for $\b \in B$. If $r_j\b=\b$ we would have
$\a_k\in B$ but all elements of $B$ except $\a_i$ are orthogonal
to $\a_i$.  This means $\a_k$ is not orthogonal to $\a_j$ and we
have $j\sim k$, $j\sim i$, and $i\sim k$ a contraction as there
are no triangles in the Dynkin diagram. We conclude that the braid
relations hold if either $\t_i$ or $\t_j$ annihilates $x_B$.

We now consider the cases in which $i\not\sim j$ with neither
$\a_i$ nor $\a_j$ being in $B$.  We wish to show
$\t_i\t_j=\t_j\t_i$.

We suppose first that both $\a_i$ and $\a_j$ are in $B^\perp$.
This means that $\t_ix_B=x_Bh_{B,i}$ and that
$\t_jx_B=x_Bh_{B,j}$.  We need only ensure that $h_{B,i}$ and
$h_{B,j}$ commute, which is Proposition \ref{prop:hBi}(ii).

Suppose now that $\a_i$ is in $B^\perp$ and $\a_j$ is not in
$B^\perp$.  In this case $\t_j\t_ix_B=\t_jx_Bh_{B,i}$. Also
$\t_jx_B=x_{r_jB}-\delta mx_B$ where $\delta$ is $0$ or $1$.  This
gives
$$\t_j\t_ix_B=(x_{r_jB}-\delta mx_B)h_{B,i}.$$
We also get $\t_i\t_jx_B=\t_ix_{r_jB}-\delta m\t_ix_B$. Notice
$\a_i\in B^\perp$ and $i\not\sim j$ imply $\a_i\in (r_jB)^\perp$.
In particular
$$\t_i\t_jx_B=x_{r_jB}h_{r_jB,i}-\delta m x_Bh_{B,i}.$$
In order for this to be $\t_j\t_ix_B$ we need
$h_{r_jB,i}=h_{B,i}$, which is satisfied by Proposition
\ref{prop:hBi}(iv).

We are left with the case in which neither $\a_i$ nor $\a_j$ is in
$B$ or in $B^\perp$.  In this case the relevant actions are $\t_i$
on $x_B$ and $\t_j$ on $x_B$.  If $r_iB=r_jB$ it is clear $\tau_i$
and $\t_j$ commute.  This gives the table
\medskip
\begin{center}
\begin{tabular}{|r|r|c|c|}
\hline
$\tau_i$ on $x_{B}$&$\tau_j$ on $x_B$&$\tau_i\tau_j x_B=\t_j\t_ix_B$\\
\hline
lower\ \ &lower\ \ & $x_{r_ir_jB}$  \\
lower\ \ &raise\ \ &$x_{r_ir_jB}-mx_{r_iB}$\\
raise\ \ &raise\ \ &$x_{r_ir_jB}-mx_{r_iB}-mx_{r_jB}+m^2x_B$  \\
\hline
\end{tabular}
\end{center}

Notice that $\a_i \not \in (r_jB)^\perp$ as $\a_i \not \in B^\perp$.
Similarly $\a_j\not \in (r_iB)^\perp$.

Suppose first that $\t_i$ and $\t_j$ both lower $B$.  By
Proposition \ref{Prop:AdmOrd}(iii) this means $\t_i$ also lowers
$r_jB$ and $\t_j $ lowers $r_iB$. Now
$$\t_i\t_jx_B=\t_ix_{r_iB}=x_{r_ir_jB}.$$
The same result occurs in the reverse order as $r_i$ and $r_j$
commute.

Suppose next that $\t_i$ and $\t_j$ both raise $B$.  Then by Lemma
\ref{Cor:AdmOrdProps}(ii), $\t_i$ raises $r_jB$ and $\t_j$ raises
$r_iB$.  In particular we have
$$\t_j\t_ix_B=\t_j(x_{r_iB}-mx_B)=x_{r_jr_iB}-mx_{r_iB}-mx_{r_jB}+m^2x_B.$$
The same is true for the reverse order.

Suppose then $\t_i$ lowers $B$ and $\t_j$ raises $B$.  By Lemma
\ref{Cor:AdmOrdProps}(i), applied to $\{\t_iB<B<\t_jB\}$, the
reflection $r_i$ also lowers $r_jB$ and $r_j$ raises $r_iB$.  This
means
$$\t_i\t_jx_B=\t_i(x_{r_jB}-mx_B)=x_{r_ir_jB}-mx_{r_iB}.$$
In the other order
$$\t_j\t_ix_B=\t_jx_{r_iB}=x_{r_jr_iB}-mx_{r_iB}.$$ These are the same.
Notice here the assumptions imply $r_iB\neq r_jB$ and $r_ir_jB\neq
B$. We conclude that $\t_i$ and $\t_j$ commute whenever $i\not\sim
j$.

\medskip
We now suppose $i\sim j$ and wish to show
$\t_i\t_j\t_i=\t_j\t_i\t_j$. Suppose first $\a_i$ and $\a_j$ are
in $B^\perp$.  Then $\t_ix_B=x_Bh_{B,i}$ and $\t_jx_B=x_Bh_{B,j}$.
The condition needed is
$h_{B,i}h_{B,j}h_{B,i}=h_{B,j}h_{B,i}h_{B,j}$, which is
Proposition \ref{prop:hBi}(iii).

Suppose now $i\sim j$ and $\a_i\in B^\perp$ but $\a_j\not \in B^{\perp}$.
We are still asssuming neither $\a_i$ nor $\a_j$ is in $B$.  The
relevant data here are the actions of $r_j$ on $B$ and $r_i$ on
${r_jB}.$ The table below handles the cases where $r_j$ lowers $B$
and $r_i$ lowers $r_jB$ as well as those where $r_j$ raises $B$
and $r_i$ raises $r_jB$.  The other cases, of $r_i$ raising $r_jB$
when $r_j$ lowers $B$ and of $r_i$ lowering $r_jB$ when $r_i$
raise $B$, are ruled out by Condition (ii) of Proposition
\ref{Prop:AdmOrd}.
\medskip
\begin{center}
\begin{tabular}{|r|r|c|c|}
\hline
$r_i$ on ${r_jB}$&$r_j$ on $B$&$\tau_i\tau_j\tau_i x_B$\\
\hline
lower\ \ &lower\ \ & $x_{r_ir_jB}h_{r_ir_jB,j}=x_{r_ir_jB}h_{B,i}$  \\
raise\ \ &raise\ \ &$x_{r_ir_jB}h_{r_ir_jB,j}-mx_B-x_{r_jB}h_{B,i}+m^2x_Bh_{B,i}$  \\
\hline
\end{tabular}
\end{center}

Notice that $\a_i\not \in (r_jB)^\perp$ as if $(\a_j,\b)\neq 0$,
then $(\a_i,r_j\b)=(\a_i,\b-(\a_j,\b) \a_j)=(\a_j,\b)\neq 0$.

Suppose first $r_j$ lowers $B$ and $r_i$ lowers $r_jB$ as in the
first row.  Then
$$\t_j\t_i\t_jx_B=\t_j\t_ix_{r_jB}=\t_jx_{r_ir_jB}=x_{r_ir_jB}h_{r_ir_jB,j}.$$
Note here $\a_j\in (r_ir_jB)^\perp$ by application of $r_ir_j$ to
$\a_i\in B^\perp$.  Also
$$\t_i\t_j\t_ix_B=\t_i\t_jx_Bh_{B,i}=\t_ix_{r_jB}h_{B,i}=x_{r_ir_jB}h_{B,i}.$$
Now the braid relation is satisfied according to Proposition
\ref{prop:hBi}(v).

Suppose $r_j$ raises $B$ and $r_i$ raises $r_jB$.
\begin{eqnarray*}
\t_i\t_j\t_ix_B & = & \t_i\t_jx_Bh_{B,i}    \\
              &=&\t_i(x_{r_jB}-mx_B)h_{B,i} \\
            &=& (x_{r_ir_jB}-mx_{r_jB}-mx_Bh_{B,i})h_{B,i}  \\
        & = & x_{r_ir_jB}h_{B,i}-mx_{r_jB}h_{B,i}-mx_Bh_{B,i}^2 \\
            &=&x_{r_ir_jB}h_{B,i}-mx_{r_jB}h_{B,i} -mx_B+m^2x_Bh_{B,i}
\end{eqnarray*}
Here we used $h_{B,i}^2=1-mh_{B,i}$. In the other order we have
\begin{eqnarray*}
\t_j\t_i\t_jx_B &=&\t_j\t_i(x_{r_jB}-mx_B)  \\
                &=&\t_j(x_{r_ir_jB}-mx_{r_jB}-mx_Bh_{B,i})  \\
                   &=&x_{r_ir_jB}h_{r_ir_jB,j}-mx_B-m(x_{r_jB}-mx_B)h_{B,i}  \\
                  &=&x_{r_ir_jB}h_{r_ir_jB,j} - mx_B -mx_{r_jB}h_{B,i} + m^2x_Bh_{B,i}.
\end{eqnarray*}
Once again we need $h_{r_ir_jB,j}=h_{B,i}$ which is Proposition
\ref{prop:hBi}(v).

\medskip
We can finally consider the case in which $i\sim j$ and neither
$\a_i$ nor $\a_j$ is in $B^\perp\cup B$.  Here relevant data are
the actions of $r_i$ and $r_j$ on $B$, where for the first row we
assume $\a_i+\a_j\not\in B$ (for otherwise, each side equals
zero).

\medskip
\begin{center}
\begin{tabular}{|r|r|c|c|}
\hline
$r_i$ on $B$&$r_j$ on $B$&$\tau_i\tau_j\tau_i x_B$\\
\hline
lower\ \ &lower\ \ & $x_{r_ir_jr_iB}$  \\
lower\ \ &raise\ \ &  done below   \\
raise\ \ &raise\ \ &$   x_{r_ir_jr_iB}-m(x_{r_jr_iB}+x_{r_ir_jB})$     \\
         &         &$\ \ \ +m^2(x_{r_jB}+x_{r_iB}) -(m^3+m)x_B$    \\
\hline
\end{tabular}
\end{center}

We start with the first row in which both $r_i$ and $r_j$ lower
$B$.  We may assume $r_iB\neq r_jB$ or $\tau_i$ and $\tau_j$ act
on $x_B$ and $x_{r_iB}$ in the same way. By Proposition
\ref{Prop:AdmOrd}(iv) and Lemma \ref{heighttworoot} all the actions we encounter are lowering
actions. Therefore,
$$\t_i\t_j\t_ix_B=\t_i\t_jx_{r_iB}=\t_ix_{r_jr_iB}=x_{r_ir_jr_iB}.$$
This gives the same result with the other product.

Next take the bottom row in which both $r_i$ and $r_j$ raise $B$.
By Lemma \ref{Cor:AdmOrdProps} (iii), the actions we encounter are
all raising actions.
\begin{eqnarray*}
\t_i\t_j\t_ix_B &=& \t_i\t_j(x_{r_iB}-mx_B)  \\
                 &=&\t_i(x_{r_jr_iB}-mx_{r_iB}-m(x_{r_jB}-mx_B)) \\
                    &=&x_{r_ir_jr_iB}-mx_{r_jr_iB}-mx_B   \\
                      &&\ \ \ -m(x_{r_ir_j B}-mx_{r_jB})+m^2(x_{r_iB}-mx_B)  \\
                        &=&x_{r_ir_jr_iB}-m(x_{r_jr_iB}+x_{r_ir_jB})   \\
                           & &\ \ \ +m^2(x_{r_jB}+x_{r_iB}) -(m^3+m)x_B.\\
\end{eqnarray*}
This also gives the same result with the other product.

We now tackle the remaining cases. Here $r_i$ lowers $B$ and $r_j$
raises $B$. There are two cases depending on how $r_j$ acts on
$r_iB$.

\medskip
\begin{center}
\begin{tabular}{|r|r|c|c|}
\hline
$r_i$ on $B$&$r_j$ on $B$&$r_j $ on $r_iB$ &  $\t_i\t_j\t_ix_B=\t_j\t_i\t_jx_B$      \\
\hline
lower\ \ &raise\ \ & raise &  $x_{r_ir_jr_iB}-mx_{r_jr_iB} -m(x_B-mx_{r_iB})$       \\
lower\ \ &raise\ \ & lower  & $x_{r_jr_ir_jB}-mx_{r_jr_iB}$
\\
\hline
\end{tabular}
\end{center}

Consider first the second row, where $r_j$ lowers $r_iB$.  By the
Lemma \ref{Cor:AdmOrdProps}(iv) applied to $r_ir_jB$, this means
$r_i$ raises $r_jr_iB$ and the remaining raising and lowering
actions can be determined by this. Notice $r_jr_iB\neq B$, for
otherwise $r_iB=r_jB$ which is not consistent with the assumption.
\begin{eqnarray*}
\t_i\t_j\t_ix_B  &=&  \t_i\t_jx_{r_iB}  \\
                   &=&\t_ix_{r_jr_iB} \\
                      &=&x_{r_ir_jr_iB}-mx_{r_jr_iB}
\end{eqnarray*}
For the other product
\begin{eqnarray*}
\t_j\t_i\t_jx_B &=& \t_j\t_i(x_{r_jB}-mx_B)  \\
                   &=&\t_j(x_{r_ir_jB}-mx_{r_iB})  \\
                      &=&x_{r_jr_ir_jB}-mx_{r_jr_iB}
\end{eqnarray*}
These are the same as indicated in the table.

For the first row suppose next that $r_j$ raises $r_iB$.  By Lemma
\ref{Cor:AdmOrdProps}(iii) applied to $r_iB$, $r_j$ raises $r_iB$,
$r_i$ raises $r_jr_iB$, $r_j$ lowers $r_ir_jB$ and $r_i$ raises
$r_jB$.  Again we use $B\ne r_jr_iB$.
\begin{eqnarray*}
\t_i\t_j\t_ix_B &=& \t_i\t_jx_{r_iB}  \\
                &=& \t_i(x_{r_jr_iB}-mx_{r_iB})  \\
                 &=&x_{r_ir_jr_iB}-mx_{r_jr_iB} -m(x_B-mx_{r_iB}).
\end{eqnarray*}

For the other product
\begin{eqnarray*}
\t_j\t_i\t_jx_B  &=& \t_j\t_i(x_{r_jB}-mx_B)  \\
                   &=&\t_j(x_{r_ir_jB}-mx_{r_jB} -mx_{r_iB})  \\
                    &=&x_{r_jr_ir_jB}-mx_B -m(x_{r_jr_iB}-mx_{r_iB}).
\end{eqnarray*}
This gives the same for either product.

These are also the same as indicated in the table finishing the
last case.  In particular Theorem \ref{th:main} has been proven.
\end{proof}

We expect that the representations obtained for the positive
monoid $A^+$ by means of our Main Theorem \ref{th:main} will be
extendible to the full Artin group $A$. Proving this is work in
progress. For type ${\rm A}_n$, all of them are, as is clear from
the BMW algebra of that type (\cite{wenzl,CGW}).

\end{document}